\documentclass[lettersize,journal]{IEEEtran}
\usepackage{amsmath,amsfonts}
\usepackage{algorithmic}
\usepackage{algorithm}
\usepackage{array}
\usepackage[caption=false,font=normalsize,labelfont=sf,textfont=sf]{subfig}
\usepackage{textcomp}
\usepackage{stfloats}
\usepackage{url}
\usepackage{verbatim}
\usepackage{graphicx}
\usepackage{cite}
\usepackage{doi}

\newtheorem{theorem}{Theorem}

\begin{document}

\title{Fast Sampling for Linear Inverse Problems of Vectors and Tensors using Multilinear Extensions}

\author{Hao Li, Dong Liang, Zixi Zhou, Zheng Xie
  \thanks{Hao Li, Dong Liang, Zixi Zhou, Zheng Xie are with College of Sciences, National University of Defense Technology, Changsha, 410073, Hunan, China (Corresponding author: liangd2018311316@163.com)}

  \thanks{This work was supported by the National Natural Science Foundation of China (Grant No.61773020). }
}


\IEEEpubid{0000--0000/00\$00.00~\copyright~2023 IEEE}

\maketitle

\begin{abstract}
  This paper studies the problem of sampling vector and tensor signals, which is the process of choosing sites in vectors and tensors to place sensors for better recovery. A small core tensor and multiple factor matrices can be used to sparsely represent a dense higher-order tensor within a linear model. Using this linear model, one can effectively recover the whole signals from a limited number of measurements by solving  linear inverse problems (LIPs).  By providing the closed-form expressions of multilinear extensions for the frame potential of pruned matrices,  we develop an algorithm named fast Frank-Wolfe algorithm (FFW) for sampling vectors and tensors with low complexity. We provide the approximation factor of our proposed algorithm for    the factor matrices that are non-orthogonal and have elements of the same sign in each row.
  Moreover, we conduct experiments to verify the higher performance and lower complexity of our proposed algorithm   for general factor matrix.  Finally, we demonstrate that  sampling by FFW and reconstruction by least squares methods  yield better results for image data compared to convCNP completion with random sampling.
\end{abstract}


\section{Introduction}
\IEEEPARstart{R}{eal} world discrete signals like audio, images and videos can be mathematically described as vectors, matrices, or higher-order tensors that reside in various modes \cite{F_3,F_4}. A target signal can be retrieved from incomplete samples by solving a linear inverse problem (LIP), which is generally considered to be linear in the remaining factor mode(s) with low-dimensional parameters \cite{F_11}.
It is crucial for applications where signal acquisition is expensive or time-consuming  to intelligently select partial entries of a signal for better reconstruction through sampling for LIPs \cite{F_14}, such as sensor signals in wireless communication \cite{F_15}, labels on medical images \cite{F_16}, and ratings in recommendation systems \cite{F_17}. Essentially a combinatorial problem, the  sampling problem is very challenging to solve, even for small-scale issues. Therefore, generating a suboptimal sampling strategy with a quality guarantee is crucial. The commonly used quality guarantee is the approximation factor, which is the ratio of the value of the suboptimal  solution to the value of the best solution.

Suppose that the samples are corrupted by independently and identically distributed (i.i.d.) noise, the quantity and quality of observed samples will have a significant impact on the reconstructed mean square error (MSE) of the unbiased least-squares (LS) solution \cite{P_21}. As a result, extensive literature has examined sampling methods to reduce reconstruction MSE while working within a limited sampling budget.

Sampling vectors has been applied in sensor placement to deploy sensors for monitoring a physical field.
The sensor placement problem was initially described as a nonconvex optimization problem by \cite{P_25}, and was relaxed into a convex problem that can be solved using interior point methods with polynomial complexity.
However, when the sensor budget was very small, this relaxation strategy worked badly. After that, greedy algorithms were created to choose sensor locations one by one by solving various local optimization problems to improve MSE performance \cite{P_21,P_22,Linear}. \cite{Linear} introduced FrameSense, which preserves the sub-modularity property and has theoretically bounded performance. The worst-case MSE function was used as the objective to select samples utilizing developed efficient algorithms in \cite{P_21}. However, the computation of the eigenspace corresponding to the chosen sensors was necessary for each greedy search, and  it is expensive when the eigenspace dimension is large. A fast MSE pursuit algorithm was recently proposed by \cite{P_22} to greedily reduce an approximate MSE criterion. However, this algorithm was still plagued by repeated matrix inverse calculations.

Although there are numerous sampling methods for vector signals, the high cost of computation and storage prevents them from being used for higher-order tensor signals. Therefore, \cite{C} developed a multi-domain frame potential-based sampling strategy to address these issues by designing samples in each order and then combining their intersection tensor as the chosen data. The proposed sampling methods utilized factor matrices directly since their  sampling matrix has the Kronecker product structure. Additionally, \cite{P_1,P_29} described three structured sub-Nyquist sampling methods for tensor signal querying. These methods were referred to as slab sampling, fiber sampling, and entry sampling, respectively. These sample techniques cannot reach arbitrary sampling budgets, despite their minimal complexity.  A fast MSE-based and unstructured sampling technique for linear-model signals ranging from vectors to tensors was recently proposed by \cite{P}. This method works well for vector signals, but it has a high computational complexity for high-order tensors.

\IEEEpubidadjcol

We make three contributions. First, we provide the closed-form expressions of multilinear extensions for the frame potential of   pruned  matrices. Second, we propose a fast algorithm (FFW) for sampling vector signals, and its solution is guaranteed to be close to the optimal one in a special case (Theorem \ref{theo2}). Third, we extend FFW to  sampling tensor signals, and also give the guarantee of solution quality for a special class of factor matrices (Theorem \ref{theo3}).
Moreover, our experiments verify that FFW is effective for more general cases where factor matrices do  not satisfy the conditions of Theorem \ref{theo2} and  Theorem \ref{theo3}.
We  experimentally compare FFW with FrameSense \cite{Linear}, Greedy FP \cite{C} and the recently proposed FMBS \cite{P} to verify the effectiveness of our proposed algorithm.

\section{Problem statements}
Throughout this paper, we use $\left( \cdot \right) ^{\dagger}$, $\left( \cdot \right) ^{T}$ and  $\left< \cdot ,\cdot \right> $ to represent Moore-Penrose pseudoinverse, transposition and the inner product, respectively. Matrix (vector) symbols are represented by upper (lower) case boldface letters, such as $\mathbf{X} (\boldsymbol{x})$.
The expectation operator is denoted by $\mathbb{E} \left\{ \cdot \right\} $. $\otimes$  represents the Kronecker product, and the main properties  used in this paper are $(\mathbf{A}\otimes\mathbf{B})(\mathbf{C}\otimes\mathbf{D})=\mathbf{A}\mathbf{C}\otimes\mathbf{B}\mathbf{D}$ and $(\mathbf{A}\otimes\mathbf{B})^\dagger=\mathbf{A}^\dagger\otimes\mathbf{B}^\dagger$.
Some important properties of Kronecker product can be seen in \cite{C_11}. 

A tensor $\mathcal{F} \in \mathbb{R}^{N_1 \times \cdots \times N_R}$ of order $R$ can be viewed as a discretized multidomain signal, with each of its entries indexed over $R$ different domains.
Two tensors $\mathcal{F} \in \mathbb{R}^{N_1 \times \cdots \times N_R}$ and $\mathcal{G} \in \mathbb{R}^{K_1 \times \cdots \times K_R}$ can be related by a multilinear system of equations as
\begin{equation}\label{eq4}
  \mathcal{F}=\mathcal{G} \bullet_1 \mathbf{U}_1 \bullet_2 \cdots \bullet_R \mathbf{U}_R,
\end{equation}
where $\left\{\mathbf{U}_r \in \mathbb{R}^{N_r \times K_r}\right\}_{r=1}^R$ represents a set of factor matrices that relates each domain of $\mathcal{F}$ and $\mathcal{G}$, and $\bullet_r$ represents the $r$-th mode product between a tensor and a matrix \cite{C_12}. Vectorizing \eqref{eq4}, we have
\begin{equation}\label{eq5}
  \boldsymbol{f}=\left(\mathbf{U}_1 \otimes \cdots \otimes \mathbf{U}_R\right) \boldsymbol{g}
\end{equation}
with $\boldsymbol{f}=\operatorname{vec}(\mathcal{F}) \in \mathbb{R}^{\widetilde{N}},\ \widetilde{N}=\prod_{r=1}^R N_r$, and $\boldsymbol{g}=\operatorname{vec}(\mathcal{G}) \in$ $\mathbb{R}^{\widetilde{K}},\ \widetilde{K}=\prod_{r=1}^R K_r$.

In this paper, we are concerned with  sampling a tensor  $\mathcal{F} $, which is equivalent to selecting entries of $\boldsymbol{f}=\operatorname{vec}({\mathcal{F}})$.   Suppose that the set of  matrices $\left\{\mathbf{U}_r\right\}_{r=1}^R$ are perfectly known, and that each of them is tall, i.e., $N_r>K_r$ for $r=1, \cdots, R$, and has full column rank.

Let $\mathcal{N}_{r}$ be the set of all row indices of the matrix $\mathbf{U}_{r}$  and $\mathcal{L}_{r}$ be the set of selected row indices from $\mathcal{N}_{r}$, $r=1,\cdots, R$.
In order to circumvent the curse of dimensionality, Ortiz-Jimenez et al. \cite{C}  defined a sampling matrix
$$\boldsymbol{\Phi}(\mathcal{L}):=\boldsymbol{\Phi}_{1}\left(\mathcal{L}_{1}\right) \otimes \cdots \otimes \boldsymbol{\Phi}_{R}\left(\mathcal{L}_{R}\right),$$
where  $\mathcal{L}=\bigcup_{r=1}^{R} \mathcal{L}_{r}$ and $\boldsymbol{\Phi}_{r}\left(\mathcal{L}_{r}\right)$ is a selection matrix for $\mathbf{U}_r$, $r=1,\cdots, R$. The notation used here refers to \cite{C}. Then, sampling a tensor  can be performed independently for each domain, that is
\begin{equation}\label{eq6}
  \begin{aligned}
    \boldsymbol{v} & =\mathbf{\Phi }(\mathcal{L} )\boldsymbol{f}
    \\&=\left( \mathbf{\Phi }_1\left( \mathcal{L} _1 \right) \otimes \cdots \otimes \mathbf{\Phi }_R\left( \mathcal{L} _R \right) \right) \left( \mathbf{U}_1\otimes \cdots \otimes \mathbf{U}_R \right) \boldsymbol{g}
    \\&=\left( \mathbf{\Phi }_1\left( \mathcal{L} _1 \right) \mathbf{U}_1\otimes \cdots \otimes \mathbf{\Phi }_R\left( \mathcal{L} _R \right) \mathbf{U}_R \right) \boldsymbol{g}.
  \end{aligned}
\end{equation}

Let $\left|\mathcal{L}_{r}\right|=L_{r}$ be the number of selected sensors per domain and $|\mathcal{L}|=\sum_{r=1}^{R} L_{r}=L$  be the total number of selected sensors.

Denote $
  \mathbf{\Psi }(\mathcal{L} )=\mathbf{\Psi }_1\left( \mathcal{L} _1 \right) \otimes \cdots \otimes \mathbf{\Psi }_R\left( \mathcal{L} _R \right)
$, in which $\mathbf{\Psi}_r\left(\mathcal{L}_r\right)=\boldsymbol{\Phi}_r\left(\mathcal{L}_r\right) \mathbf{U}_r$, $r=1,2,\cdots,R$. In the rest of this paper,  we will omit the $\left( \mathcal{N} \right)$ and $\left( \mathcal{N} _r \right)$ of $ \mathbf{\Psi }(\mathcal{N} )$ and  $\mathbf{\Psi}_r\left(\mathcal{N}_r\right)$ when it is clear from the context.

Notice that \eqref{eq6} is overdetermined, by least squares, we can estimate the core as
$$
  \hat{\boldsymbol{g}}=\mathbf{\Psi }^{\dagger}(\mathcal{L} )\boldsymbol{v}=\left[ \left( \mathbf{\Psi }_1\left( \mathcal{L} _1 \right) \right) ^{\dagger}\otimes \cdots \otimes \left( \mathbf{\Psi }_R\left( \mathcal{L} _R \right) \right) ^{\dagger} \right] \boldsymbol{v},
$$
and then reconstruct $\hat{\boldsymbol{f}}$ by \eqref{eq5}.



If the measurements collected in $\boldsymbol{v}$ are perturbed by zero-mean white Gaussian noise with unit variance, then the least-squares solution has  the Fisher information matrix given by $\mathbf{T}(\mathcal{L})=\mathbb{E}\left\{(\boldsymbol{g}-\hat{\boldsymbol{g}})(\boldsymbol{g}-\hat{\boldsymbol{g}})^{T}\right\}=\mathbf{\Psi}^{T}(\mathcal{L})\mathbf{\Psi}(\mathcal{L})$, which determines the quality of the estimators $\hat{\boldsymbol{g}}$.
Thus, a sampling problem can be posed as a discrete optimization problem that finds the best sampling subset by optimizing a scalar function of $\mathbf{T}(\mathcal{L})$. The  most popular scalar function is the mean squared error $\operatorname{MSE}(\mathbf{\Psi}(\mathcal{L}))=\operatorname{tr}\left\{\mathbf{T}^{-1}(\mathcal{L})\right\}$, which is  difficult to minimize as it is neither convex, nor submodular.

In order to design an efficient sampling strategy, we introduced the frame potential (FP) \cite{C_25} of the matrix $\mathbf{\Psi}(\mathcal{L})$ as the cost function to replace MSE, which  is defined as $\operatorname{FP}(\mathbf{\Psi}(\mathcal{L})):=\operatorname{tr}\left\{\mathbf{T}^T(\mathcal{L}) \mathbf{T}(\mathcal{L})\right\}$. Following the Lemma 2 in \cite{Linear}, the MSE is bounded by the FP as,
$$
  c_1 \frac{\operatorname{FP}(\mathbf{\Psi}(\mathcal{L}))}{\lambda_{\max }^2\{\mathbf{T}(\mathcal{L})\}} \leq \operatorname{MSE}(\mathbf{\Psi}(\mathcal{L})) \leq c_2 \frac{\operatorname{FP}(\mathbf{\Psi}(\mathcal{L}))}{\lambda_{\min }^2\{\mathbf{T}(\mathcal{L})\}},
$$
where $c_1$, and $c_2$ are constants that depend on the data model.
From the above bound, it is clear that  one can minimize $\operatorname{MSE}(\mathbf{\Psi}(\mathcal{L}))$ by minimizing $\operatorname{FP}(\mathbf{\Psi}(\mathcal{L}))$.


Following \cite{C}, the frame potential of $\mathbf{\Psi}(\mathcal{L})$ can be expressed as
\begin{equation}\label{eq7}
  \mathrm{FP(}\mathbf{\Psi }(\mathcal{L} ))=\prod_{r=1}^R{\mathrm{FP}}\left( \mathbf{\Psi }_r\left( \mathcal{L} _r \right) \right)=
  \prod_{r=1}^R{\mathrm{tr}}\left\{ \mathbf{T}_{r}^{T}\left( \mathcal{L} _r \right) \mathbf{T}_r\left( \mathcal{L} _r \right) \right\}
\end{equation}
where $\mathbf{T}_r\left(\mathcal{L}_r\right) =\mathbf{\Psi}_r^T\left(\mathcal{L}_r\right)  \mathbf{\Psi}_r\left(\mathcal{L}_r\right) $.

For brevity, let
$$
  F\left( \mathcal{L} \right) =\prod_{r=1}^R{F_r\left( \mathcal{L} _r \right)}:=\prod_{r=1}^R{\mathrm{FP}\left( \mathbf{\Psi }_r\left( \mathcal{L} _r \right) \right)},
$$
where $\mathcal{L}_{i} \cap\mathcal{L}_{j}=\emptyset $ for $i \neq j$.


Denote $N=\sum_{r=1}^R{N_r}$ and $K=\sum_{r=1}^R{K_r}$. The sampling problem can be formulated as
\begin{equation}
  \begin{aligned}\label{p1}
              & \underset{\mathcal{L} _1,\cdots ,\mathcal{L} _R}{\min}F\left( \mathcal{L} \right) =\prod_{r=1}^R{\mathrm{FP}}\left( \mathbf{\Psi }_r\left( \mathcal{L} _r \right) \right) \\
    s.t.\quad & \sum_{r=1}^R{\left| \mathcal{L} _r \right|}=L,                                                                                                                            \\
              & \left| \mathcal{L} _r \right|\ge K_r, r=1,2,\cdots ,R,                                                                                                                    \\
  \end{aligned}
\end{equation}
where $L$ is a positive integer satisfying $L \geq K$, which represents the number of sensors selected.
Sampling vectors and tensors is Problem \eqref{p1} with $R=1$ and $R\geq 2$, respectively.

Let $
  \widetilde{L}=\prod_{r=1}^R{L_r}
$
represents the total number of samples collected using the Kronecker-structured sampler mentioned above. It should be noted that the cardinality constraint in Problem \ref{p1} limits the total number of selected sensors to $L= \sum_{r=1}^R{L_r}$.
If the total number of selected sensors is limited by  $\widetilde{L}$, the complexity of  the   near optimal solvers \cite{c15,  C16}  of Problem \ref{p1} is $\mathcal{O}(N^5)$, which will not be extended to large-scale problems.


\section{Preliminaries}
In this section, due to the difficulties in analysis techniques for discrete functions, we give the multilinear extension of $F(\mathcal{L} )$, which transform $F(\mathcal{L} )$ into a continuous function.
The set function $F(\mathcal{L} )$ is defined on the vertices of the hypercube $[0,1]^{N}$ and each set $\mathcal{L}$  is a vertex. The multilinear extension $\widetilde{F}$ of $F(\mathcal{L} )$ is defined as
$$
  \widetilde{F}(\boldsymbol{x}):=\underset{\mathcal{L} \sim \boldsymbol{x}}{\mathbb{E}}\left\{ F(\mathcal{L} ) \right\} =\sum_{\mathcal{L} \subseteq \mathcal{N}}{F(\mathcal{L} )\prod_{j\in \mathcal{L}}{x_j}\prod_{j\notin \mathcal{L}}{1}-x_j},
$$
where $\mathcal{L} \sim \boldsymbol{x}$ represents that $\mathcal{L}$ is drawn from the independent distribution with marginals $\boldsymbol{x} \in [0,1]^{N}$.  One can see that for any set $\mathcal{L} $ and its indicator vector $\mathbf{1}_\mathcal{L} $, $\widetilde{F}\left(\mathbf{1}_\mathcal{L}\right)=F(\mathcal{L} )$. The following theorem shows the closed-form expression of $\widetilde{F}(\boldsymbol{x}) $.

\begin{theorem}\label{theo0}
  Let $p^r_{ab}$ be the $(a,b)$-element of $\mathbf{\Psi}_r$. The multilinear extension of $F(\mathcal{L} )$ can be expressed  as
  $$
    \widetilde{F}(\boldsymbol{x})=\prod_{r=1}^R{\widetilde{F}_r\left( \boldsymbol{x}^r \right)}=\prod_{r=1}^R{\left\{ \sum_{i=1}^{K_r}{\sum_{j=1}^{K_r}{\left[ \sum_{n=1}^{N_r}{x_{n}^{r}p_{ni}^{r}}p_{nj}^{r} \right] ^2}} \right\}},
  $$
  where $
    \boldsymbol{x}=\left( \boldsymbol{x}^1,\boldsymbol{x}^2,\cdots ,\boldsymbol{x}^R \right)^T$ and $\boldsymbol{x}^r=\left( x_{1}^{r},x_{2}^{r},\cdots ,x_{N_r}^{r} \right)$, $r=1,\cdots,R$.
\end{theorem}
\begin{IEEEproof}
  Note that $\widetilde{F}(\boldsymbol{x})$ is the expected value of $F(\mathcal{L} )$ over the independent distribution with marginals $\boldsymbol{x}$, we have
  $$
    \begin{aligned}
       & \quad \widetilde{F}(\boldsymbol{x})                                                                                                                                                                                                                            \\
       & =\sum_{\mathcal{L} \subseteq \mathcal{N}}{\left\{ \prod_{r=1}^R{F_r\left( \mathcal{L} _r \right)}\prod_{r=1}^R{\prod_{j\in \mathcal{L} _r}{x_{j}^{r}}}\prod_{r=1}^R{\prod_{j\in \mathcal{N} _r\backslash\mathcal{L} _r}{\left( 1-x_{j}^{r} \right)}} \right\}}
      \\
       & =\sum_{\mathcal{L} \subseteq \mathcal{N}}{\left\{ \prod_{r=1}^R{\left[ F_r\left( \mathcal{L} _r \right) \prod_{j\in \mathcal{L} _r}{x_{j}^{r}}\prod_{j\in \mathcal{N} _r\backslash\mathcal{L} _r}{\left( 1-x_{j}^{r} \right)} \right]} \right\}}
      \\
       & =\prod_{r=1}^R{\left\{ \sum_{\mathcal{L} _r\subseteq \mathcal{N} _r}{F_r\left( \mathcal{L} _r \right) \prod_{j\in \mathcal{L} _r}{x_{j}^{r}}\prod_{j\in \mathcal{N} _r\backslash\mathcal{L} _r}{\left( 1-x_{j}^{r} \right)}} \right\}}
      \\
       & =\prod_{r=1}^R{\widetilde{F_r}\left( \boldsymbol{x}^r \right)}.
    \end{aligned}
  $$

  By the definitions of $\widetilde{F}_r(\boldsymbol{x}^r)$ and $F(\mathcal{L} )$, we obtain
  $$
    \begin{aligned}
      \widetilde{F}_r(\boldsymbol{x}^r) & =\sum_{\mathcal{L} _r\subseteq \mathcal{N} _r}{F_r(\mathcal{L} _r)\prod_{j\in \mathcal{L} _r}{x_{j}^{r}}\prod_{j\in \mathcal{N} _r\backslash \mathcal{L} _r}{1}-x_{j}^{r}}                                                                                                    \\
                                        & =\sum_{\mathcal{L} _r\subseteq \mathcal{N} _r}{\mathrm{tr}\left\{ \mathbf{T}_{r}^{T}\left( \mathcal{L} _r \right) \mathbf{T}_r\left( \mathcal{L} _r \right) \right\} \prod_{j\in \mathcal{L} _r}{x_{j}^{r}}\prod_{j\in \mathcal{N} _r\backslash \mathcal{L} _r}{1}-x_{j}^{r}} \\
                                        & =\sum_{a=1}^{K_r}{\sum_{b=1}^{K_r}{\sum_{\mathcal{L} _r\subseteq \mathcal{N} _r}{\left[ \mathbf{T}_r\left( \mathcal{L} _r \right) _{ab} \right] ^2\prod_{j\in \mathcal{L} _r}{x_{j}^{r}}\prod_{j\in \mathcal{N} _r\backslash \mathcal{L} _r}{1}-x_{j}^{r}}}}                  \\
                                        & =\sum_{a=1}^{K_r}{\sum_{b=1}^{K_r}{\underset{\mathcal{L} _r\sim \boldsymbol{x}^r}{\mathbb{E}}\left\{ \mathbf{T}_r\left( \mathcal{L} _r \right) _{ab} \right\} ^2}},                                                                                                           \\
    \end{aligned}
  $$
  where $
    \mathbf{T}_r\left( \mathcal{L} _r \right) _{ab}
  $ denotes the $(a,b)$-element of $
    \mathbf{T}_r\left( \mathcal{L} _r \right)$. Then, we have
  $$
    \begin{aligned}
      \underset{\mathcal{L} _r\sim \boldsymbol{x}^r}{\mathbb{E}}\left\{ \mathbf{T}_r\left( \mathcal{L} _r \right) _{ab} \right\} ^2 & =\left[ \underset{\mathcal{L} _r\sim \boldsymbol{x}^r}{\mathbb{E}}\left\{ \mathbf{T}_r\left( \mathcal{L} _r \right) _{ab} \right\} \right] ^2  \\
                                                                                                                                    & =\left[ \underset{\mathcal{L} _r\sim \boldsymbol{x}^r}{\mathbb{E}}\left\{ \sum_{n\in \mathcal{L} _r}{p_{na}^{r}p_{nb}^{r}} \right\} \right] ^2 \\
                                                                                                                                    & =\left[ \sum_{n=1}^{N_r}{x_{n}^{r}p_{na}^{r}}p_{nb}^{r} \right] ^2.                                                                            \\
    \end{aligned}
  $$
  Therefore, we conclude that
  $$
    \widetilde{F}_r(\boldsymbol{x}^r)=\sum_{i=1}^{K_r}{\sum_{j=1}^{K_r}{\left[ \sum_{n=1}^{N_r}{x_{n}^{r}p_{ni}^{r}}p_{nj}^{r} \right] ^2}}.
  $$
  This completes the proof of Theorem \ref{theo0}.
\end{IEEEproof}

According to Theorem \ref{theo1},  Problem \ref{p1} is transformed into
$$
  \begin{aligned}
              & \underset{\boldsymbol{x}^1,\boldsymbol{x}^2,\cdots ,\boldsymbol{x}^R}{\min}\widetilde{F}(\boldsymbol{x}) \\
    s.t.\quad & \left\| \boldsymbol{x} \right\| _1=L,                                                                    \\
              & \left\| \boldsymbol{x}^r \right\| _1\ge K_r,\  r=1,2,\cdots ,R,                                          \\
              & \boldsymbol{x}\in \left\{ 0,1 \right\} ^N,                                                               \\
  \end{aligned}
$$
whose objective function is derivable.
The general technique for minimizing $\widetilde{F}(\boldsymbol{x})$ is the continuous greedy algorithm \cite{M_50}, which is a slight modification of the Frank-Wolfe algorithm (FW) \cite{M_13}, with a fixed step size $\delta=1 / N^2$. In each iteration, the algorithm takes a step $\boldsymbol{x}^{(t+1)}=\boldsymbol{x}^{(t)}+\delta \boldsymbol{h}^{(t)}$ in the direction $$
  \boldsymbol{h}^{(t)}=\mathrm{arg}\min _{\boldsymbol{h}^{\prime}\in \mathcal{P}}\left. \langle \boldsymbol{h}^{\prime},\nabla \widetilde{F}(\boldsymbol{x}^{\left( t \right)}) \right. \rangle,
$$
where $\mathcal{P}$ denotes the polytope corresponding to the family of feasible sets.
However, for any step size $\delta=1 / t$ with $t \geq 2$, this algorithm terminates in $\mathcal{O}\left(t\right)$ iterations and gives a fractional solution. In order to round the fractional solution to  an integral solution, the pipage rounding technique \cite{M_50} is needed, which requires $\mathcal{O}\left(N^2\right)$ function calls.
Obviously,  this algorithm with $t\geq 2$  is poorly suited for large-scale problems because of  the excessive complexity of pipage rounding.

Here, motivated by the continuous greedy algorithm, we design a  fast algorithm with low complexity which can be seen as $t=1$ and  $
  \boldsymbol{x}^{\left( 0 \right)}=\varepsilon \mathbf{1}_N\ (\varepsilon>0)
$ in FW. Moreover, through the
closed-form expression of  $\widetilde{F}(\boldsymbol{x})$, we further reduce the complexity of the algorithm and can directly obtain an integral solution, which corresponds to a unique set.


\section{Algorithm to  Sampling Vectors}

In this  section, we focus on   vector signals ($R=1$) and propose the Fast Frank-Wolfe algorithm (FFW) for computing suboptimal solution of Problem \eqref{p1}.

Since there is only one factor matrix for vector signals, we have  $\mathbf{\Psi}=\mathbf{\Psi}_1,\ N=N_1$ and $ K=K_1$.
For convenience, denote  by $
  \boldsymbol{x}=\left( x_1,x_2,\cdots x_N \right)^T
$, $p_{ab}$ the  $(a,b)$-element of $\mathbf{\Psi}$ and  $\boldsymbol{p}_i$ the $i$-th column of $\mathbf{\Psi}$. The multilinear extension of $F(\mathcal{L} )$ can be expressed as $$
  \widetilde{F}(\boldsymbol{x})=\sum_{i=1}^K{\sum_{j=1}^K{\left[ \sum_{n=1}^N{x_n p_{ni}p_{nj}} \right] ^2}}.
$$
According to this expression, we give the following algorithm to sample vectors.

\subsection{The Algorithm}
Note that for any $\varepsilon >0$, the  partial  derivative of $ \widetilde{F}(\boldsymbol{x})$ at $ \varepsilon \mathbf{1}_N$ are
$$
  \begin{aligned}
    \left. \frac{\partial \widetilde{F}_r(\boldsymbol{x})}{\partial x_t} \right|_{\boldsymbol{x}=\varepsilon \mathbf{1}_N} & =2\varepsilon \sum_{i=1}^K{\sum_{j=1}^K{\left[ \sum_{n=1}^N{p_{ni}p_{nj}} \right]}}p_{ti}p_{tj}
    \\
                                                                                                                           & =2\varepsilon \sum_{i=1}^K{\sum_{j=1}^K{\left( {\boldsymbol{p}_i}^T\boldsymbol{p}_j \right) p_{ti}p_{tj}}},\ t=1,\cdots ,N.
  \end{aligned}
$$
The main idea of our algorithm is to select the smallest $L$ partial derivatives among $N$ partial derivatives at $ \varepsilon \mathbf{1}_N$, and take the corresponding rows as the selected rows.

\begin{algorithm}[H]
  {\bf {Input:}} $\mathbf{\Psi}=(p_{ab})\in \mathbb{R}^{N \times K}$; $L$\;

  1: Compute $
    m_{ij } =\boldsymbol{p}_i^T\boldsymbol{p}_j,\ i=1,\cdots ,K,\ j=1,\cdots ,K
  $\;

  2: Compute $
    d_n =\sum_{i=1}^K{\sum_{j=1}^K{ m_{ij }  p_{ni}p_{nj}}},\ n=1,\cdots ,N
  $\;

  3: Select the $L$ elements with the smallest values from $
    \left\{ d_1,d_ 2  ,\cdots ,d_ N  \right\}
  $, denoted as $
    \left\{ d_ {l_1}  ,d_ {l_2 } ,\cdots ,d_ {l_{L}}\right\}
  $\;

  4: {\bf {Return:}} $
    \mathcal{L} =\left\{ l_1,l_2,\cdots ,l_{L} \right\}$
  \caption{Fast Frank-Wolfe Algorithm for Vector Signals.}
  \label{alg1}
\end{algorithm}



The computational complexity of Steps 1 and 2 is $\mathcal{O} \left( NK^2 \right)$ and that of Step 3 is $\mathcal{O} \left( N \right)$. Therefore, the total computational complexity of Algorithm \ref{alg1} is $\mathcal{O} \left( NK^2\right)$. We compare the computational complexity with the following sampling methods for vector signals: FrameSense \cite{Linear}, minimum nonzero eigenvalue pursuit (MNEP) \cite{P_21}, maximal projection on minimum eigenspace (MPME) \cite{P_21}, fast MSE pursuit-based sampling (fastMSE) \cite{P_22}, and Fast MSE-Based Sampling (FMBS) \cite{P}. The computational complexities of those methods and the proposed FFW are illustrated in Table \ref{table1}.
One can see that FFW has  the lowest complexity compared with other methods. Moreover,  the complexity of FFW does not depend on $L$. Therefore, the advantage of FFW is more prominent  when $K \ll  L$. In these methods, only FrameSense and FFW have theoretically bounded performance.

\begin{table}[h]

  \centering
  \caption{Comparison of computational complexity of different sampling methods}
  \begin{tabular}{cccc}
    \hline \hline Method & FrameSense \cite{Linear}          & MNEP \cite{P_21}                  & MPME \cite{P_21}                  \\
    Complexity           & $\mathcal{O}\left(N^3\right)$     & $\mathcal{O}\left(N L K^3\right)$ & $\mathcal{O}\left(N L K^2\right)$ \\
    \hline Method        & fastMSE \cite{P_22}               & FMBS  \cite{P}                    & FFW                               \\
    Complexity           & $\mathcal{O}\left(N L K^2\right)$ & $\mathcal{O}\left(N L^2\right)$   & $\mathcal{O}\left(N K^2\right)$   \\
    \hline \hline
  \end{tabular}\label{table1}
\end{table}


\subsection{Near-optimality of FFW for Vector Signals}
Here, we give  the approximation factor of Algorithm \ref{alg1}.

\begin{theorem}\label{theo1}
  Denote by ${\mathcal{L}^*}$ an optimal solution of Problem \eqref{p1}. If $\mathbf{\Psi}$ is a non-orthogonal matrix and the elements of each row in $\mathbf{\Psi}$ have the same sign,
  Algorithm \ref{alg1} returns a set $\mathcal{L}^{\prime}$ such that $
    G(\mathcal{S} ^{\prime})>\frac{N}{N+L}G(\mathcal{S} ^*)
  $, where $
    \mathcal{S} ^{\prime}=\mathcal{N} /\mathcal{L} ^{\prime},\ \mathcal{S} ^*=\mathcal{N} /\mathcal{L} ^*
  $ and
  $G\left( \mathcal{S} \right) =F\left( \mathcal{N} \right) -F\left( \mathcal{N} /\mathcal{S} \right) $.

\end{theorem}

\begin{IEEEproof}
  Let  $\mathcal{H} =\left\{ \boldsymbol{z}\in \left\{ 0,1 \right\} ^N|\left\| \boldsymbol{z} \right\| _1=N-L \right\}
  $. Denote by $\widetilde{G}(\boldsymbol{z})$ the multilinear extension of $G\left( \mathcal{S} \right)$ and $\boldsymbol{z}\in\mathcal{H} $.

  Note that $
    \widetilde{G}\left( \mathbf{0}_N \right) =G\left( \emptyset \right) =0
  $   and $\widetilde{G}\left( \boldsymbol{z} \right)$ is a quadratic function, $\widetilde{G}\left( \boldsymbol{z} \right)$ can be express as
  \begin{equation} \label{eq8}
    \widetilde{G}(\boldsymbol{z})=\boldsymbol{d} ^T\boldsymbol{z}+\frac{1}{2}\boldsymbol{z}^T \mathbf{H} \boldsymbol{z},
  \end{equation}
  where $\boldsymbol{d}=(d_1,d_2,\cdots,d_N)^T$,
  \begin{equation} \label{eq9}
    d_t=\left. \frac{\partial \widetilde{G}\left( \boldsymbol{z} \right)}{\partial z_t} \right|_{\boldsymbol{z}=\mathbf{0}_N}=2\sum_{i=1}^K{\sum_{j=1}^K{\left( \boldsymbol{p}_{i}^{T}\boldsymbol{p}_j\cdot p_{ti}p_{tj} \right)}}
  \end{equation}
  and $\mathbf{H}=\{h_{s,t}\}_{N\times N}$,
  \begin{equation}\label{eq10}
    h_{s,t}\left. =\frac{\partial ^2\widetilde{G}\left( \boldsymbol{z} \right)}{\partial z_t\partial z_s} \right|_{\boldsymbol{z}=\mathbf{0}_N}=-2\sum_{i=1}^K{\sum_{j=1}^K{\left( p_{si}p_{sj}\cdot p_{ti}p_{tj} \right)}}.
  \end{equation}
  Substituting \eqref{eq9} and \eqref{eq10} in \eqref{eq8}, we obtain
  $$
    \begin{aligned}
      \widetilde{G}(\boldsymbol{z}) & =2\sum_{i=1}^K{\sum_{j=1}^K{\left( \boldsymbol{p}_{i}^{T}\boldsymbol{p}_j\cdot \sum_{n=1}^N{z_np_{ni}p_{nj}} \right)}}
      \\
                                    & \quad -\sum_{i=1}^K{\sum_{j=1}^K{\left[ \sum_{n_1=1}^N{\sum_{n_2=1}^N{\left( z_{n_1}p_{n_1i}p_{n_1j}\cdot z_{n_2}p_{n_2i}p_{n_2j} \right)}} \right]}}
      \\
                                    & =\sum_{i=1}^K{\sum_{j=1}^K{\left[ 2\boldsymbol{p}_{i}^{T}\boldsymbol{p}_j\sum_{n=1}^N{z_np_{ni}p_{nj}}-\left( \sum_{n=1}^N{z_np_{ni}p_{nj}} \right) ^2 \right]}}.
    \end{aligned}
  $$
  Let $
    M=\underset{\boldsymbol{z}\in \mathcal{H}}{\max}\,\,\boldsymbol{d}^T\boldsymbol{z}
  $.
  One can see that $
    \boldsymbol{d}^T\mathbf{1}_{\mathcal{S} ^{\prime}}=\boldsymbol{d}^T\left( \mathbf{1}_N-\mathbf{1}_{\mathcal{L} ^{\prime}} \right) =M
  $     as $\mathcal{L}^{\prime}$ is obtained by Algorithm \ref{alg1}.
  Let $\alpha _{ij}=\boldsymbol{p}_{i}^T\boldsymbol{p}_j$ and $y_{ij}=\sum_{n=1}^N{z_np_{ni}p_{nj}}$, $i,j=1,\cdots,K$. Then $\widetilde{G}(\boldsymbol{z})$ can be express as
  $$
    \widehat{G}(\boldsymbol{y})=\sum_{i=1}^K{\sum_{j=1}^K{\left[ 2\alpha _{ij}\cdot y_{ij}-{y_{ij}}^2 \right]}}.
  $$
  where $
    \boldsymbol{y}=\left( y_{11},\cdots ,y_{1K,}\cdots \cdots ,y_{K1},\cdots ,y_{KK} \right) ^T$.

  Note that $\boldsymbol{z}\in \left\{ 0,1 \right\} ^N$ and the elements of each row in $\mathbf{\Psi}$ have the same sign, we have $\alpha _{ij} \geq  y_{ij} \geq  0$. Since $\mathbf{\Psi}$ is a non-orthogonal matrix, we get $\sum_{a=1}^K{\sum_{b=1}^K{{\alpha _{ab}}^2}}>0$.

  We first derive the biggest value that $\widehat{G}\left( \boldsymbol{y} \right)$ can obtain, denoted by $G_{\max}$, that is, to solve the following optimization problem,
  $$
    \begin{aligned}
       & \max_{\boldsymbol{y}} \,\,\,\widehat{G}(\boldsymbol{y})                              \\
       & s.t.\quad 2\sum_{i=1}^K{\sum_{j=1}^K{\left( \alpha _{ij}\cdot y_{ij} \right)}}\leq M \\
       & \qquad \,\, \alpha _{ij}\ge y_{ij}\ge 0, i,j=1,\cdots K.                             \\
    \end{aligned}
  $$

  By simple calculation, we get
  $$
    \begin{aligned}
      T_{\max} & =\sum_{i=1}^K{\sum_{j=1}^K{\left[ 2\alpha _{ij}\cdot y_{ij}^{*}-\left( y_{ij}^{*} \right) ^2 \right]}}
      \\
               & =M-\frac{M^2}{4\sum_{a=1}^K{\sum_{b=1}^K{{\alpha _{ab}}^2}}},
    \end{aligned}
  $$
  in which $$
    y_{ij}^{*}=\frac{M\alpha _{ij}}{2\sum_{a=1}^K{\sum_{b=1}^K{{\alpha _{ab}}^2}}},\
    i,j=1,\cdots,K.$$
  Thus,
  \begin{equation}\label{eq15}
    \widetilde{G}\left( \boldsymbol{z} \right) \le M-\frac{M^2}{4\sum_{a=1}^K{\sum_{b=1}^K{{\alpha _{ab}}^2}}}.
  \end{equation}

  Now, we give an lower bound $T_{\min}$ of $\widehat{G}(\boldsymbol{y})$ in the case $ \boldsymbol{z}=\mathbf{1}_{\mathcal{S} ^{\prime}}$.
  We can obtain $G_{\min}$ by solving the following optimization problem,
  $$
    \begin{aligned}
           & \min_{\boldsymbol{y}} \,\,\,\widehat{G}(\boldsymbol{y})
      \\
      s.t. & \quad 2\sum_{i=1}^K{\sum_{j=1}^K{\left( \alpha _{ij}\cdot y_{ij} \right)}}=M
      \\
           & \qquad \,\,\alpha _{ij}\ge y_{ij}\ge 0,i,j=1,\cdots K.
    \end{aligned}
  $$

  Sort ${\alpha_{ij}}$, $i,j=1,2,\cdots K$ and renumber ${\alpha_{ij}}$ by $
    {\alpha_1} \geq  {\alpha _2}\geq  \cdots \geq  {\alpha_{K^2}} $. There exist integer $k$  and real number $\delta$ such that $\sum_{i=k}^{K^2}{\alpha _i}^2\leq  \frac{M}{2}<\sum_{i=k-1}^{K^2}{\alpha_i}^2$ and $
    \sum_{i=k}^{K^2}{{\alpha _i}^2}+\alpha _{k-1}\delta =\frac{M}{2}
  $, where $0\leq \delta<\alpha_{k-1}$. Then we have
  $$
    \begin{aligned}
      G_{\min} & =\left( 2\alpha _{k-1}\cdot \delta -\delta ^2 \right) +\sum_{i=k}^K{\left( 2\alpha _i\cdot \alpha _i-{\alpha _i}^2 \right)}
      \\
               & \quad \ge \alpha _{k-1}\cdot \delta +\sum_{i=k}^K{{\alpha _i}^2}=\frac{M}{2}
    \end{aligned}
  $$
  and \begin{equation}\label{eq16}
    \widetilde{G}\left( \mathbf{1}_{\mathcal{S} ^{\prime}} \right) \ge \frac{M}{2}.
  \end{equation}

  Note that $
    \frac{M}{\sum_{i=1}^K{\sum_{j=1}^K{{\alpha _{ij}}^2}}}=\frac{\underset{\boldsymbol{x} \in \mathcal{H}}{\max}\,\,D^T\boldsymbol{x}}{\,\,\frac{1}{2}D^T\mathbf{1}_N}\geq \frac{2\left( N-L \right)}{N}
  $, combining with \eqref{eq15} and \eqref{eq16}, and we have
  $$
    \frac{G(\mathcal{S} ^{\prime})}{G(\mathcal{S} ^*)}=\frac{\widetilde{G}\left( \mathbf{1}_{\mathcal{S} ^{\prime}} \right)}{\widetilde{G}\left( \mathbf{1}_{\mathcal{S} ^*} \right)}>\frac{2}{4-\frac{M}{\sum_{i=1}^K{\sum_{j=1}^K{{\alpha _{ij}}^2}}}}\ge \frac{N}{N+L}.
  $$

  This completes the proof of Theorem \ref{theo1}.

\end{IEEEproof}

Note that the approximation factors of continuous greedy algorithm and FrameSense with regard to $G\left( \mathcal{S} \right)$ are both $1-\frac{1}{e}$. Thus, if $L<\frac{N}{e-1}$, the approximation factor of FFW given by Theorem \ref{theo1} is better than continuous greedy algorithm and FrameSense.

According to the proof of Theorem 2 in \cite{Linear}, we directly derive a  bound with regard to the frame potential from Theorem \ref{theo1}.
\begin{theorem} \label{theo2}
  Suppose that $\mathbf{\Psi}\in \mathbb{R}^{N \times K}$  satisfies the condition of Theorem \ref{theo1}. Denote by $\mathcal{L}^*$ an optimal solution of Problem \eqref{p1} with $R=1$. The set $\mathcal{L}^{\prime}$ obtained from Algorithm \ref{alg1}  is near optimal with respect to FP as $
    F\left( \mathcal{L}^{\prime} \right) \leq  \gamma F\left( \mathcal{L} ^* \right)
  $ with $
    \gamma =\frac{1}{N+L}\left( F\left(\mathcal{N}  \right) \frac{KL}{L_{\min}^{2}}+N \right)
  $, and $
    L_{\min}=\min_{\left| \mathcal{L} \right|=L} \sum_{i\in \mathcal{L}}{\left\| u_i \right\| _{2}^{2}}
  $, being $u_i$ the $i$-th row of $\mathbf{\Psi}$.
\end{theorem}

It should be noted that the condition of Theorem \ref{theo1} is not necessary for the  performance of FFW, that is, even if the factor matrix $\mathbf{\Psi}$ does not satisfy the conditions of Theorem \ref{theo1}, FFW can still perform  well in terms of MSE. We verify this conclusion through experiments, refer to the third experiment in  VI (A).


\section{Algorithm to Sampling Tensors}

In this section, we focus on the   general situation of Problem \eqref{p1} with $R \geq  2$. For convenience, let $p_{ab}^r$ be the $(a,b)$-element of $\mathbf{\Psi}_r$ and $\boldsymbol{p}_i^r$ be the $i$-th column of $\mathbf{\Psi}_r$, $r=1,\cdots,R$.

According to the closed-form expression of $\widetilde{F}(\boldsymbol{x}) $ given by Theorem \ref{theo1}, we give the following algorithm to sample tensors.
\subsection{The Algorithm}

For any $\varepsilon >0$, we have
$$
  \begin{aligned}
     & \quad \frac{\varepsilon}{2}\frac{1}{\widetilde{F}(\varepsilon \mathbf{1}_N)}\left. \cdot \frac{\partial \widetilde{F}(\boldsymbol{x})}{\partial x_{t}^{r}} \right|_{\boldsymbol{x}=\varepsilon \mathbf{1}_N}                                                                                            \\
     & =\frac{\varepsilon}{2}\frac{\prod_{a\ne r}{\widetilde{F_a}\left( \varepsilon \mathbf{1}_{N_a} \right)}}{\widetilde{F}(\varepsilon \mathbf{1}_N)}\cdot \left. \frac{\partial \widetilde{F_r}\left( \boldsymbol{x}^r \right)}{\partial x_{t}^{r}} \right|_{\boldsymbol{x}^r=\varepsilon \mathbf{1}_{N_r}} \\
     & =\frac{\varepsilon}{2}\frac{1}{\widetilde{F_r}\left( \varepsilon \mathbf{1}_{N_r} \right)}\cdot \left. \frac{\partial \widetilde{F_r}\left( \boldsymbol{x}^r \right)}{\partial x_{t}^{r}} \right|_{\boldsymbol{x}^r=\varepsilon \mathbf{1}_{N_r}}                                                       \\
     & =\frac{1}{\sum_{i=1}^{K_r}{\sum_{j=1}^{K_r}{\left( {\boldsymbol{p}_{i}^{r}}^T\boldsymbol{p}_{j}^{r} \right) ^2}}}\cdot \sum_{i=1}^{K_r}{\sum_{j=1}^{K_r}{\left( {\boldsymbol{p}_{i}^{r}}^T\boldsymbol{p}_{j}^{r} \right)}}p_{ti}^{r}p_{tj}^{r}                                                          \\
  \end{aligned}
$$

Analogous to  Algorithm \ref{alg1}, the main idea of Algorithm \ref{alg2} is also to select the smallest $L$ partial derivatives among $N$ partial derivatives at $ \varepsilon \mathbf{1}_N$.

\begin{algorithm}[H]
  {\bf {Input:}}$\mathbf{\Psi}_r=(p^r_{ab})\in \mathbb{R}^{N_r \times K_r},\ r=1,\cdots,R$; $L$

  1: {\bf {For $r=1\dots R$ do}}

  2:\qquad Compute $
    m^r_{ij } ={\boldsymbol{p}^r_i}^T \boldsymbol{p}^r_j,\ i,j=1,\cdots ,K_r
  $\;

  3:\qquad Compute $
    F^r=\sum_{i=1}^{K_r}{\sum_{j=1}^{K_r}{\left( m_{ij}^{r} \right) ^2}}
  $\;

  4:\qquad Compute $
    d_{n}^{r}=\frac{1}{F^r}\sum_{i=1}^{K_r}{\sum_{j=1}^{K_r}{m_{ij}^{r}p_{ni}^{r}p_{nj}^{r}}},$

  \qquad\quad $n= 1,\cdots ,N_r
  $, and denote $\mathcal{D}_r =
    \left\{ d_{1}^{r},d_{2}^{r},\cdots ,d_{N_r}^{r}\right\}$\;



  5: {\bf {End}}

  6: Select the $L$ elements with the smallest values $\bigcup_{r=1}^R{\left\{ d_{l_{1}^{r}}^{r},\cdots ,d_{l_{L_r}^{r}}^{r} \right\}}$ from $\bigcup_{r=1}^R{\mathcal{D} _r}$ such that  $\sum_{r=1}^R{L_r}=L$ and $L_r\geq K_r$, $r=1,\cdots,R$\;

  7: {\bf {Return:}} $\mathcal{L} =\bigcup_{r=1}^R\left\{ l_{1}^{r},\cdots ,l_{L_r}^{r} \right\}$ \;


  \caption{ Fast Frank-Wolfe Algorithm for Tensor Signals}
  \label{alg2}
\end{algorithm}

The computational complexity of Step 1 to Step 5 is $\mathcal{O} \left( N_{\max}K_{\max}^2 \right)$ with $
  N_{\max}=\underset{r}{\max}\,N_r$ and $K_{\max}=\underset{r}{\max}\,K_r
$. The computational complexity of Step 6 is $\mathcal{O} \left(  N\right)$. Therefore, the total computational complexity of Algorithm \ref{alg2} is $\mathcal{O} \left( N_{\max}K_{\max}^2 \right)$, which is much lower than   Greedy FP $\left(\mathcal{O} \left( N_{\max}^2K_{\max} \right)\right)$  and  UB-FMBS  $
  \left( \mathcal{O} \left( \prod_{r=1}^R{N_r} \right) \right)
$ proposed in \cite{C} and \cite{P}, respectively.

\subsection{Near-optimality of FFW for Tensor Signals}

Theorem \ref{theo3} gives the  approximation
factor of  for Algorithm \ref{alg2} for a special class of factor matrices. Moreover, the Algorithm \ref{alg2} performs well in more general cases, which can refer to the third experiment of VI (A).


\begin{theorem}\label{theo3}
  Denote by ${\mathcal{L}^*}$ an optimal solution of Problem \eqref{p1}.   If for any $r=1,\cdots,R$, $\mathbf{\Psi}_r$ is a non-orthogonal matrix
  and the elements of each row in $\mathbf{\Psi}_r$ have the same sign, then the set $\mathcal{L}^{\prime}=\bigcup_{r=1}^{R} \mathcal{L}^{\prime}_{r}$ obtained from Algorithm \ref{alg2}  is near optimal with respect to FP as $$
    F\left( \mathcal{L} ^{\prime} \right) <e^{M-\frac{M^2}{2R}+o\left( \left\| \mathbf{1}_{\mathcal{N} \backslash\mathcal{L} ^{\prime}} \right\| ^3 \right) -o\left( \left\| \mathbf{1}_{\mathcal{N} \backslash\mathcal{L} ^*} \right\| ^3 \right)}F\left( \mathcal{L} ^* \right)
  $$
  with
  $$
    M=2\sum_{r=1}^R{\left\{ \frac{\sum_{i=1}^{K_r}{\sum_{j=1}^{K_r}{\left( {\boldsymbol{p}_{i}^{r}}^T\boldsymbol{p}_{j}^{r} \right) \sum_{t\in \mathcal{S} _{r}^{\prime}}{p_{ti}^{r}p_{tj}^{r}}}}}{F_r\left( \mathcal{N} \right)} \right\} \in \left( 0,2R \right)}.
  $$

\end{theorem}
\begin{IEEEproof}
  Let $
    \mathcal{S} ^{\prime}=\bigcup_{r=1}^R{\mathcal{S} _{r}^{\prime}}$  and $\mathcal{S} ^*=\bigcup_{r=1}^R{\mathcal{S} _{r}^{*}}
  $, in which $\mathcal{S} _{r}^{\prime}=\mathcal{N} _r\backslash\mathcal{L} _{r}^{\prime}$ and $\mathcal{S} _{r}^{*}=\mathcal{N} _r\backslash\mathcal{L} _{r}^{*} $, $r=1,\cdots,R$.

  In order to represent $\widetilde{F} $ as the sum of $R$ terms, we denote
  $$
    \begin{aligned}
      W\left( \boldsymbol{z} \right) & =\log \left\{ \widetilde{F}(\mathbf{1}_N) \right\} -\log \left\{ \widetilde{F}(\mathbf{1}_N-\boldsymbol{z}) \right\}                                                                    \\
                                     & =\sum_{r=1}^R{\log \left\{ \widetilde{F_r}\left( \mathbf{1}_{N_r} \right) \right\}}-\sum_{r=1}^R{\log \left\{ \widetilde{F_r}\left( \mathbf{1}_{N_r}-\boldsymbol{z}^r \right) \right\}} \\
                                     & =\sum_{r=1}^R{\log \left\{ \sum_{i=1}^{K_r}{\sum_{j=1}^{K_r}{\left[ \sum_{n=1}^{N_r}{p_{ni}^{r}}p_{nj}^{r} \right] ^2}} \right\}}                                                       \\
                                     & \quad -\sum_{r=1}^R{\log \left\{ \sum_{i=1}^{K_r}{\sum_{j=1}^{K_r}{\left[ \sum_{n=1}^{N_r}{\left( 1-z_{n}^{r} \right) p_{ni}^{r}}p_{nj}^{r} \right] ^2}} \right\}},                     \\
    \end{aligned}
  $$
  where $
    \boldsymbol{z}=\left( \boldsymbol{z}^1,\boldsymbol{z}^2,\cdots ,\boldsymbol{z}^R \right)^T$ and $\boldsymbol{z}^r=\left( z_{1}^{r},z_{2}^{r},\cdots ,z_{N_r}^{r} \right)$, $r=1,\cdots,R$.




  For convenience,   denote

  $$
    D_{t}^{r}=\left. \frac{\partial \widetilde{F_r}\left( \mathbf{1}_{N_r}-\boldsymbol{z}^r \right)}{\partial z_{t}^{r}} \right|_{\boldsymbol{z}^r=\mathbf{0}_{N_r}}=-2\sum_{i=1}^{K_r}{\sum_{j=1}^{K_r}{\left( {\boldsymbol{p}_{i}^{r}}^T\boldsymbol{p}_{j}^{r} \right) p_{ti}^{r}p_{tj}^{r}}},
  $$
  $$
    F^r=\widetilde{F_r}\left( \mathbf{1}_{N_r} \right) =\sum_{i=1}^{K_r}{\sum_{j=1}^{K_r}{\left( {\boldsymbol{p}_{i}^{r}}^T\boldsymbol{p}_{j}^{r} \right) ^2}}
  $$
  and
  $$
    B_{st}^{r}=\frac{\partial ^2\widetilde{F_r}\left( \mathbf{1}_{N_r}-\boldsymbol{z}^r \right)}{\partial z_{t}^{r}\partial z_{s}^{r}}=2\sum_{i=1}^{K_r}{\sum_{j=1}^{K_r}{\left( p_{ti}^{r}p_{tj}^{r} \right) \left( p_{si}^{r}p_{sj}^{r} \right)}}.
  $$
  Since $\mathbf{\Psi}_r$ is a non-orthogonal matrix, $F^r>0$. Thus,
  $$
    \left. \frac{\partial W\left( \boldsymbol{z} \right)}{\partial z_{t}^{r}} \right|_{\boldsymbol{z}=\mathbf{0}_N}=-\frac{D_{t}^{r}}{F^r},
  $$
  $$
    \left.
    \frac{\partial ^2W\left( \boldsymbol{z} \right)}{\partial z_{t}^{r}\partial z_{s}^{r}} \right|_{\boldsymbol{z}=\mathbf{0}_N}=-\frac{B_{st}^{r}\cdot F^r-D_{s}^{r}D_{t}^{r}}{\left( F^r \right) ^2}.
  $$
  Using Maclaurin expansion, we get
  $$
    \begin{aligned}
      W\left( \boldsymbol{z} \right) & =-\sum_{r=1}^R{\left\{ \frac{1}{2\left( F^r \right) ^2}\sum_t^{N_r}{\sum_s^{N_r}{z_{t}^{r}z_{s}^{r}\left( B_{st}^{r}\cdot F^r-D_{s}^{r}D_{t}^{r} \right)}} \right\}} \\
                                     & \quad +\sum_{r=1}^R{\left\{ -\frac{1}{F^r}\sum_{t=1}^{N_r}{z_{t}^{r}D_{t}^{r}} \right\}}+o\left( \left\| \boldsymbol{z} \right\| ^3 \right)                          \\
                                     & =\sum_{r=1}^R{\frac{1}{F^r}\left( -\sum_{t=1}^{N_r}{z_{t}^{r}D_{t}^{r}}-\,\,\frac{1}{2}\sum_t^{N_r}{\sum_s^{N_r}{z_{t}^{r}z_{s}^{r}B_{st}^{r}}} \right)}             \\
                                     & \quad +\sum_{r=1}^R{\frac{1}{2\left( F^r \right) ^2}\left( -\sum_{t=1}^{N_r}{z_{t}^{r}D_{t}^{r}} \right) ^2}+o\left( \left\| \boldsymbol{z} \right\| ^3 \right) .    \\
    \end{aligned}
  $$

  Let $\left| \mathcal{L}^{\prime} _r \right|=L^{\prime}_r$, $\left| \mathcal{L}^* _r \right|=L^*_r$,
  $r=1,\cdots ,R
  $, where $I$ is an integer.

  Denote $$
    M_r=-\sum_{t\in \mathcal{S} _{r}^{\prime}}{D_{t}^{r}},\quad M_{r}^{*}=\underset{\left| \mathcal{S} _r \right|=N_r-L_{r}^{*}}{\max}\,\,-\sum_{\mathcal{S} _r}{D_{t}^{r}},
  $$
  and
  $$
    M=\sum_{r=1}^R{\frac{M_r}{F^r}},\quad M^*=\sum_{r=1}^R{\frac{M_{r}^{*}}{F^r}}
  $$

  Clearly, we have $M^*\leq M$ as $
    \mathcal{L} _{r}^{\prime}=\mathcal{N} _r\backslash\mathcal{S} _{r}^{\prime}
  $ is obtained by Algorithm \ref{alg2}.



  Note that $
    F^r=-\frac{1}{2}\sum_{i=1}^{N_r}{D_{i}^{r}}
  $  and  the elements of each row in $\mathbf{\Psi}_r$ have the same sign. Then  we have $0\leq M_r(M_{r}^{*})\le 2F^r$ and $0\leq M(M^*)\leq 2R$.

  We first give the upper bound of $ W\left( \mathbf{1}_{\mathcal{S}^*} \right)$. By the proof of Theorem \ref{theo1},  for $
    \boldsymbol{z}=\mathbf{1}_{\mathcal{S} ^*}
  $
  and $r=1,\cdots,R$ we have
  $$\left( -\sum_{t=1}^{N_r}{z_{t}^{r}D_{t}^{r}} \right)\leq M_r^*,$$
  $$
    \left( -\sum_{t=1}^{N_r}{z_{t}^{r}D_{t}^{r}} \right) -\frac{1}{2}\sum_t^{N_r}{\sum_s^{N_r}{z_{t}^{r}z_{s}^{r}B_{st}^{r}}}<M_{r}^{*}-\frac{\left( M_{r}^{*} \right) ^2}{4F^r}.
  $$
  Then we obtain
  \begin{equation}
    \begin{aligned}
      W\left( \mathbf{1}_{\mathcal{S}^*} \right) & <\sum_{r=1}^R{\left\{ \frac{1}{2}\left( \frac{M_{r}^{*}}{F^r} \right) ^2+\frac{1}{F^r}\left( M_{r}^{*}-\frac{\left( M_{r}^{*} \right) ^2}{4F^r} \right) \right\}}
      \\&
      \quad +o\left( \left\| \mathbf{1}_{\mathcal{S}^*} \right\| ^3 \right)
      \\&
      <\sum_{r=1}^R{\left\{ \frac{1}{4}\left( \frac{M_{r}^{*}}{F^r} \right) ^2+\frac{M_{r}^{*}}{F^r} \right\}}+o\left( \left\| \mathbf{1}_{\mathcal{S}^*} \right\| ^3 \right)
      \\&
      =\frac{1}{4}\sum_{r=1}^R{\left( \frac{M_{r}^{*}}{F^r}+2 \right) ^2}-R+o\left( \left\| \mathbf{1}_{\mathcal{S}^*} \right\| ^3 \right).
    \end{aligned}
  \end{equation}
  Denote  $Q_{max}$ the  value of
  $$
    \begin{aligned}
      \underset{y_1,y_2,\cdots y_R}{\max} & \,\,\sum_{r=1}^R{\left( y_r+2 \right) ^2}
      \\
      s.t. \quad                          & 0\le y_r\le 2
      \\
                                          & \sum_{r=1}^R{y_r}\le M.
    \end{aligned}
  $$
  Let $M=2m+\delta$, in which  $m$ is an integer and $\delta\in [0,2)$. One can verify that $Q_{max}= m\left( 2+2 \right) ^2+\left( \delta +2 \right) ^2+\left( R-m-1 \right) \cdot (0+2)^2 =12m+\delta ^2+4\delta +4R$. Thus,

  \begin{equation}\label{eq11}
    \begin{aligned}
      W\left( \mathbf{1}_{\mathcal{S}^*} \right) & < \frac{1}{4}\left( 12m+\delta ^2+4\delta +4R \right)-R+o\left( \left\| \mathbf{1}_{\mathcal{S}^*} \right\| ^3 \right)
      \\&
      <\frac{1}{4}\left( 12m+6\delta +4R\right) -R +o\left( \left\| \mathbf{1}_{\mathcal{S}^*} \right\| ^3 \right)                                                        \\&=\frac{3}{2}M+o\left( \left\| \mathbf{1}_{\mathcal{S}^*} \right\| ^3 \right).
    \end{aligned}
  \end{equation}

  Now, we derive the lower bound of $
    W\left( \mathbf{1}_{\mathcal{S} ^{\prime}} \right)
  $. Following the proof of Theorem \ref{theo1}, for $
    \boldsymbol{z}=\mathbf{1}_{\mathcal{S} ^{\prime}}
  $ and $r=1,\cdots,R$, we get
  $$
    \left( -\sum_{t=1}^{N_r}{z_{t}^{r}D_{t}^{r}} \right) =M_r,$$
  $$\quad \left( -\sum_{t=1}^{N_r}{z_{t}^{r}D_{t}^{r}} \right) -\frac{1}{2}\sum_t^{N_r}{\sum_s^{N_r}{z_{t}^{r}z_{s}^{r}B_{st}^{r}}}> \frac{M_r}{2}.
  $$
  Thus, we obtain
  \begin{equation}\label{eq12}
    \begin{aligned}
      W\left( \mathbf{1}_{\mathcal{S}^{\prime}} \right) & >\sum_{r=1}^R{\left\{\frac{1}{2\left( F^r \right) ^2}\left( M_r \right) ^2+\frac{1}{F^r}\left( \frac{M_r}{2} \right) \right\}}
      \\&\quad +o\left( \left\| \mathbf{1}_{\mathcal{S}^{\prime}}\right\| ^3 \right)
      \\&
      =\frac{1}{2}\sum_{r=1}^R{\left( \frac{M_r}{F^r}+\frac{1}{2} \right) ^2}-\frac{R}{8}+o\left( \left\| \mathbf{1}_{\mathcal{S}^{\prime}} \right\| ^3 \right).
    \end{aligned}
  \end{equation}

  Denote by $Q_{min}$ the  value of
  $$
    \begin{aligned}
      \underset{y_1,y_2,\cdots y_R}{\min} & \,\,\sum_{r=1}^R{\left( y_r+\frac{1}{2} \right) ^2}
      \\
      s.t. \quad                          & 0\le y_r\le 2
      \\
                                          & \sum_{r=1}^R{y_r}\le M.
    \end{aligned}
  $$
  Obviously, $Q_{min}=R\left( \frac{M}{R}+\frac{1}{2} \right) ^2 $. Then we have

  $$
    \begin{aligned}W\left( \mathbf{1}_{\mathcal{S}^{\prime}} \right)&>\frac{R}{2}\left( \frac{M}{R}+\frac{1}{2} \right) ^2 -\frac{R}{8}    +o\left( \left\| \mathbf{1}_{\mathcal{S}^{\prime}} \right\| ^3 \right)                                                                                                    \\&=\frac{M^2}{2R}+\frac{M}{2}+o\left( \left\| \mathbf{1}_{\mathcal{S}^{\prime}} \right\| ^3 \right).  \end{aligned}
  $$

  By \eqref{eq11} and \eqref{eq12}, we have

  $$
    \begin{aligned}
       & \quad  \log \frac{\widetilde{F}\left( \mathbf{1}_{\mathcal{L} ^*} \right)}{\widetilde{F}\left( \mathbf{1}_{\mathcal{L}^{\prime}} \right)}                                        \\ & =-\log \left\{ \widetilde{F}\left( \mathbf{1}_N-\mathbf{1}_{\mathcal{S} ^{\prime}} \right) \right\} +\log \left\{ \widetilde{F}\left( \mathbf{1}_N-\mathbf{1}_{\mathcal{S} ^*} \right) \right\}\\
       & = W\left( \mathbf{1}_{\mathcal{S} ^{\prime}} \right) -W\left( \mathbf{1}_{\mathcal{S} ^*} \right)                                                                                \\
       & >\frac{M^2}{2R}+\frac{M}{2}-\frac{3}{2}M+o\left( \left\| \mathbf{1}_{\mathcal{S} ^{\prime}} \right\| ^3 \right) -o\left( \left\| \mathbf{1}_{\mathcal{S} ^*} \right\| ^3 \right) \\
       & =\frac{M^2}{2R}-M+o\left( \left\| \mathbf{1}_{\mathcal{S} ^{\prime}} \right\| ^3 \right) -o\left( \left\| \mathbf{1}_{\mathcal{S} ^*} \right\| ^3 \right) .                      \\
    \end{aligned}
  $$

  Thus,
  $$
    \begin{aligned}
       & \quad \frac{F\left( \mathcal{L} ^* \right)}{F\left( \mathcal{L} ^{\prime} \right)}=                                                                                                                                                                                                \frac{\widetilde{F}\left( \mathbf{1}_{\mathcal{L} ^*} \right)}{\widetilde{F}\left( {\mathcal{L} ^{\prime}} \right)}>e^{\frac{M^2}{2R}-M+o\left( \left\| \mathbf{1}_{\mathcal{S} ^{\prime}} \right\| ^3 \right) -o\left( \left\| \mathbf{1}_{\mathcal{S} ^*} \right\| ^3 \right)} \\
       & \Rightarrow F\left( \mathcal{L} ^{\prime} \right)                                    <e^{M-\frac{M^2}{2R}+o\left( \left\| \mathbf{1}_{\mathcal{S} ^{\prime}} \right\| ^3 \right) -o\left( \left\| \mathbf{1}_{\mathcal{S} ^*} \right\| ^3 \right)}F\left( \mathcal{L} ^* \right) .                                                                                                                                                                                                                                                                                  \\
    \end{aligned}
  $$

  Since $
    F^r=\widetilde{F_r}\left( \mathbf{1}_N \right) =F_r\left( \mathcal{N} \right)
  $, we have
  $$
    M=\sum_{r=1}^R{\frac{M_r}{F^r}}=\sum_{r=1}^R{\left\{ \frac{-\sum_{t\in \mathcal{S} _{r}^{\prime}}{D_{t}^{r}}}{F_r\left( \mathcal{N} \right)} \right\}}.
  $$

  This completes the proof of Theorem \ref{theo3}.
\end{IEEEproof}

Different from  sampling vector signals, the algorithm of sampling  tensor signals and the analysis of theoretical boundaries are more complicated. The approximate factor given by Theorem \ref{theo2} for  sampling vectors is more accurate than the bound given by Theorem \ref{theo3} for sampling tensor, and does not contain infinitesimal terms.


\section{Numerical Results}
\subsection{Comparison with traditional sampling methods}
In this section, we experimentally compare the performance of FFW with FrameSense \cite{Linear}, Greedy FP \cite{C}, FMBS \cite{P} and Random sampling. FrameSense  is a classical method to sampling vectors, which uses greedy algorithm to optimize FP, and also gives the approximate factor. Greedy FP  extends FrameSense to sampling tensors. FMBS is the latest algorithm for sampling vectors, and its computational complexity and solution quality are the best  so far.

\begin{figure*}[!t]
  \centering
  \subfloat[Vector signals.]{\includegraphics[width=2.5in]{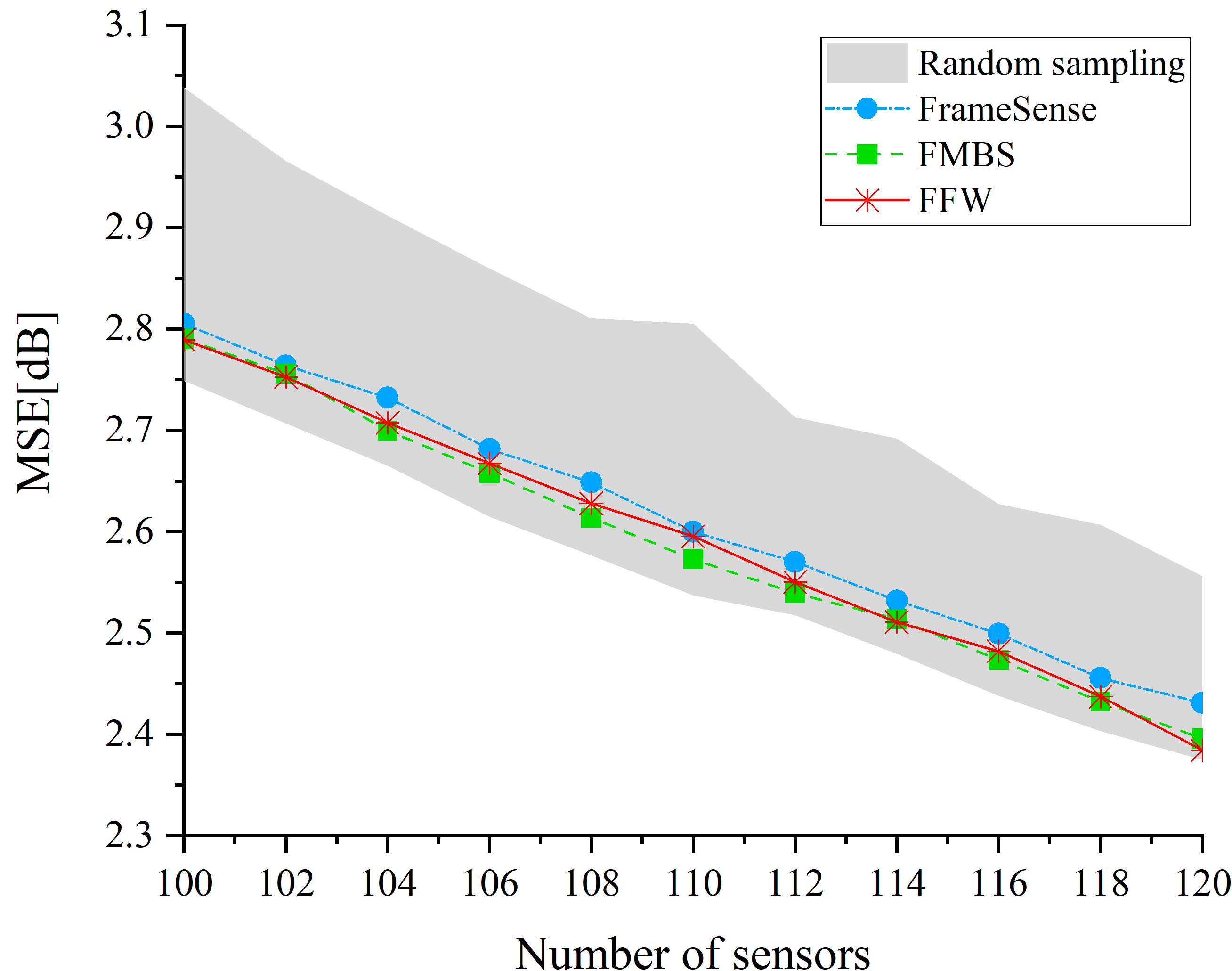}%
    \label{1_a}}
  \hfil
  \subfloat[Tensor signals.]{\includegraphics[width=2.5in]{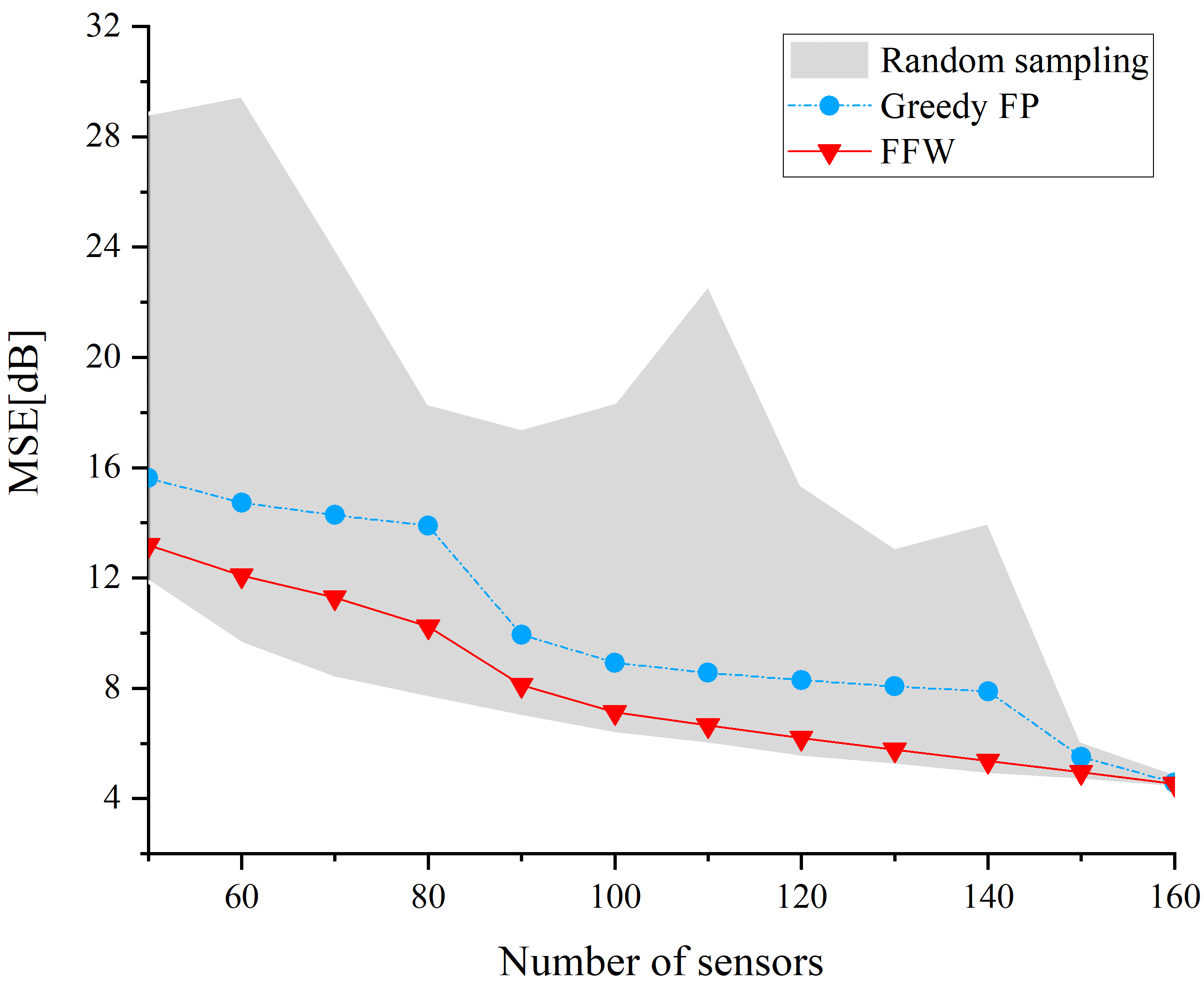}%
    \label{1_b}}
  \caption{Comparison in terms of MSE of FFW, FrameSense, FMBS, Greedy FP  and random sampling for vector and tensor signals whose factor matrices satisfy the conditions of Theorem  \ref{theo1} and Theorem \ref{theo3}.}
  \label{fig1}
\end{figure*}

In the first experiment,  for vector signals, we use 100 different instances and compute the average MSE as a function of $L$.  We randomly generate a  matrix $\mathbf{\Psi} \in \mathbb{R}^{200 \times 40}$ satisfying the conditions of Theorem \ref{theo1} for each of these instances. Then, we  compare the  performance  of  FFW with  random sampling, FrameSense  and FMBS.
Random sampling  is based on randomly selecting rows of $\mathbf{\Psi}$ for 500 times.
Since the performance is measured in terms of MSE,  the lower the curve, the higher the performance.
The results are shown in Figure \ref{1_a}. The shaded area represent the interval of random sampling.
Figure \ref{1_a} shows that FFW performs as well as FMBS, and has a slight advantage over FrameSense.


In the second experiment,  we perform the experiment with 3-order tensor signals. We also use 100 different instances to compare the performance of FFW with  random sampling  and Greedy FP. In each instance, we randomly generate three matrices $\left\{ \mathbf{\Psi }_i\in \mathbb{R} ^{N_i\times K_i} \right\} _{i=1}^{R=3}$ satisfying the conditions of Theorem \ref{theo3}  with $N_1=50, N_2=60, N_3=70, K_1=10, K_2=20, K_3=15$
and solve Problem \eqref{p1} for different number of sensors. The results are shown in Figure \ref{1_b}. It can be seen that  FFW performs much better than Greedy FP in terms of MSE. The results of FFW is closer to the optimal results of random sampling when  $L $ increases.

\begin{figure*}[!t]
  \centering
  \subfloat[Vector signals.]{\includegraphics[width=2.5in]{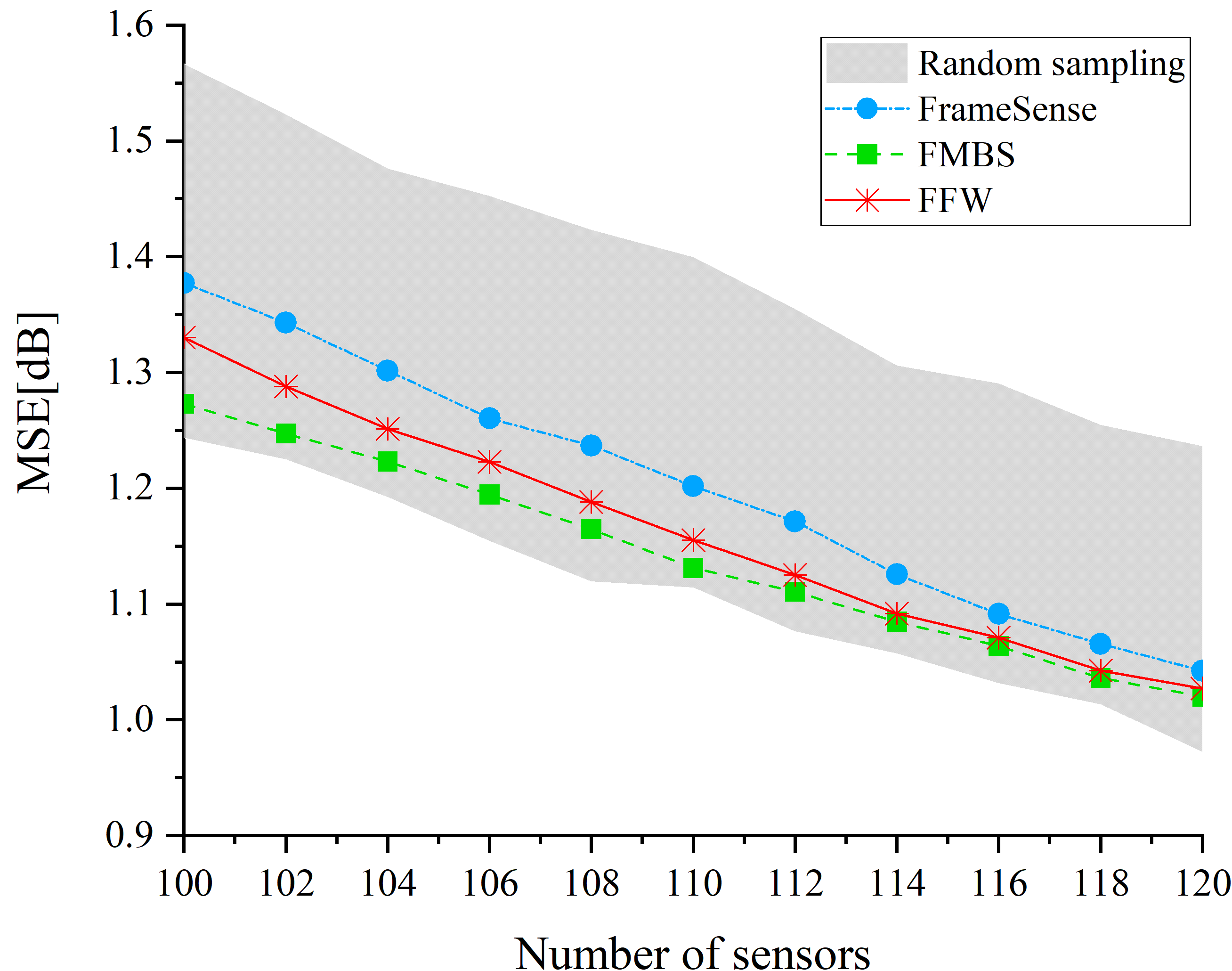}%
    \label{5_a}}
  \hfil
  \subfloat[Tensor signals.]{\includegraphics[width=2.5in]{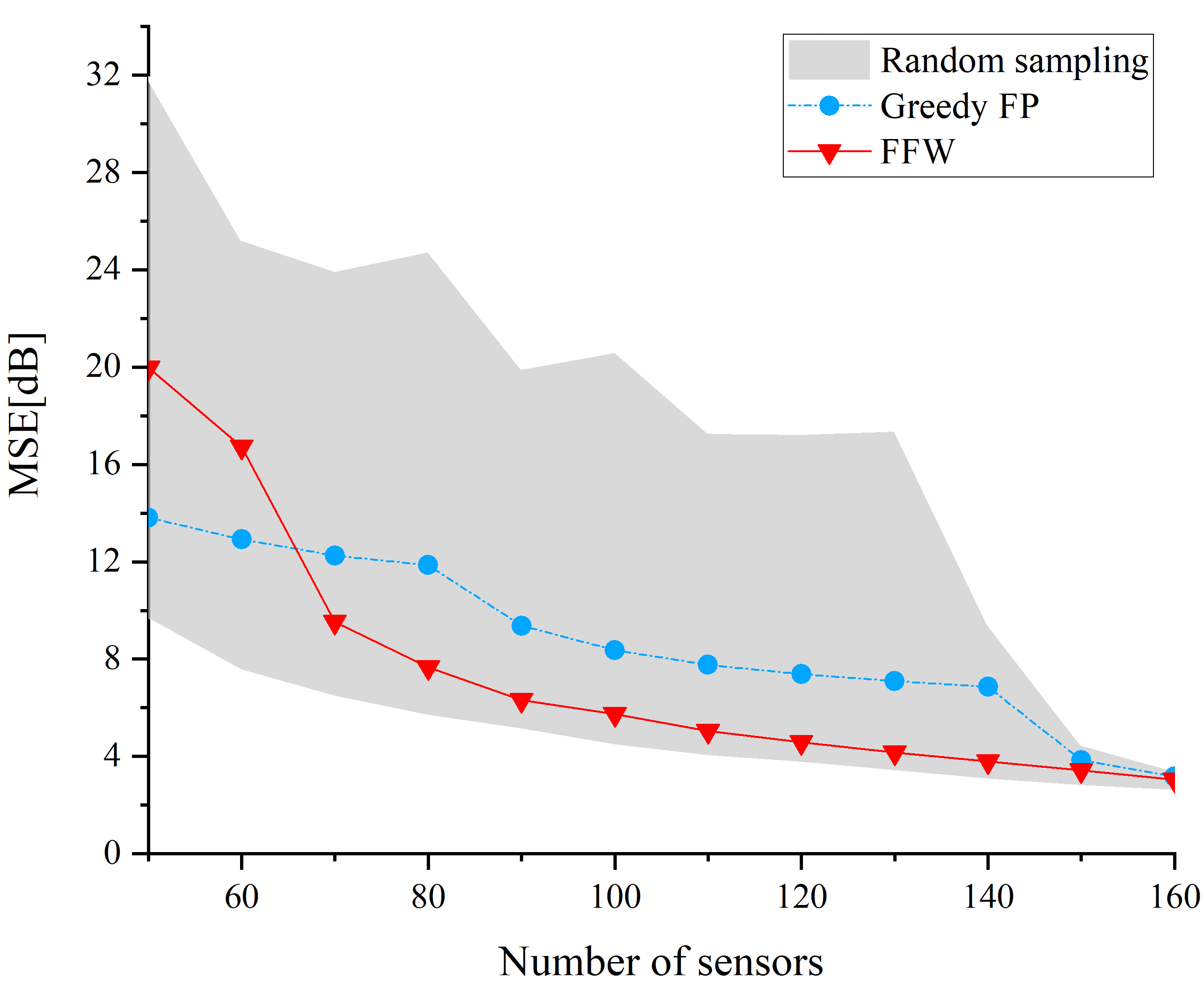}%
    \label{5_b}}
  \caption{Comparison in terms of MSE of FFW, FrameSense, FMBS, Greedy FP  and random sampling for vector and tensor signals whose factor matrices do not satisfy the conditions of Theorem  \ref{theo1} and Theorem  \ref{theo3}.}
  \label{fig5}
\end{figure*}

In the third experiment, we repeat the first and second experiments. The difference is that the matrix generated here does not satisfy the conditions of Theorem \ref{theo1} and Theorem  \ref{theo3}. The results are shown in Figure \ref{fig5}.
It can be seen that  for vector signals, the MSE of FFW is between FMBS and FrameSense. For tensor signals, FFW is better than Greedy FP when $L$ is not too small.  When $L$ is large, the advantage of FFW is more prominent for both vector and tensor signals.



\begin{figure*}[!t]
  \centering
  \subfloat[Vector signals.]{\includegraphics[width=2.5in]{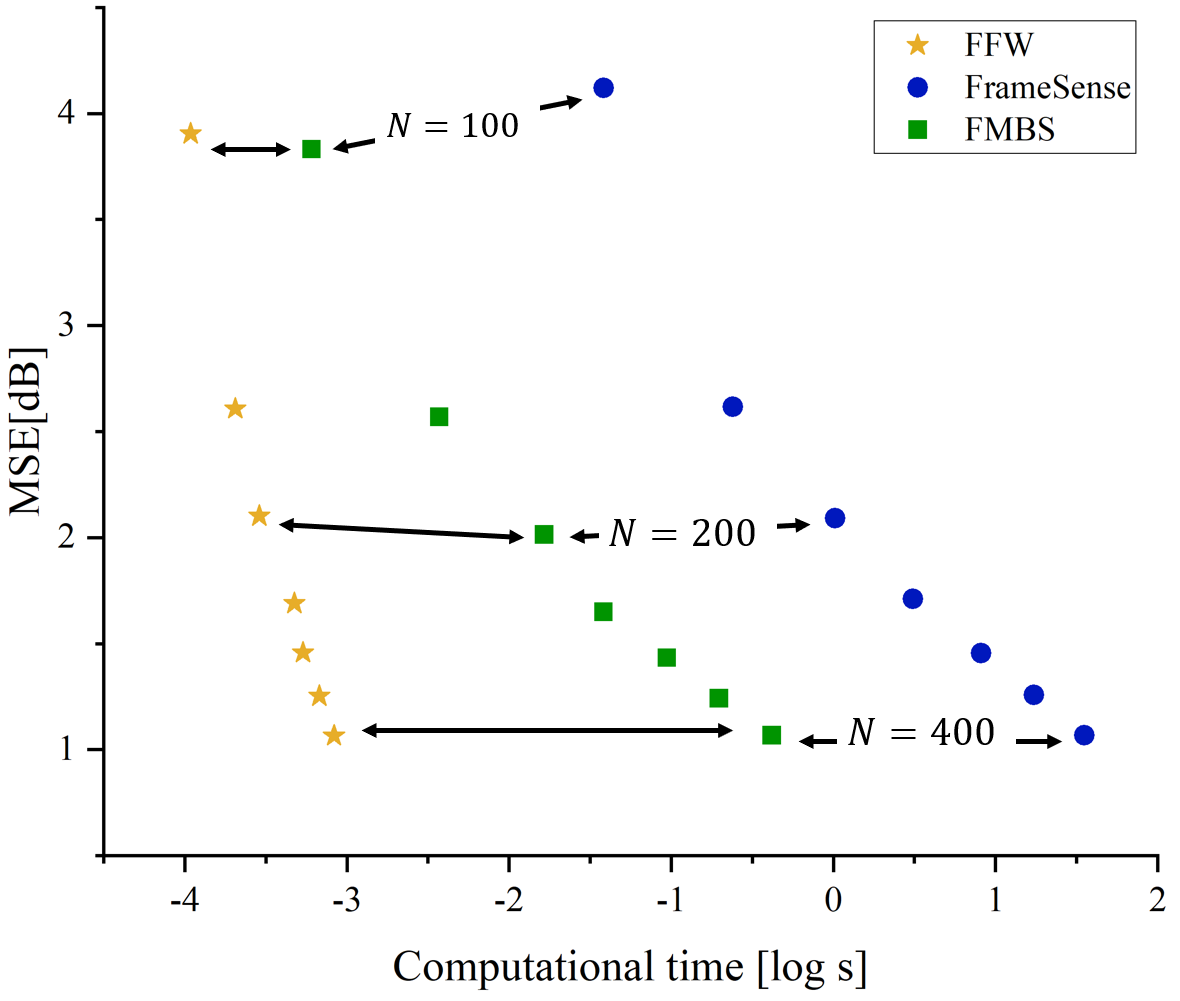}%
    \label{2_a}}
  \hfil
  \subfloat[Tensor signals.]{\includegraphics[width=2.5in]{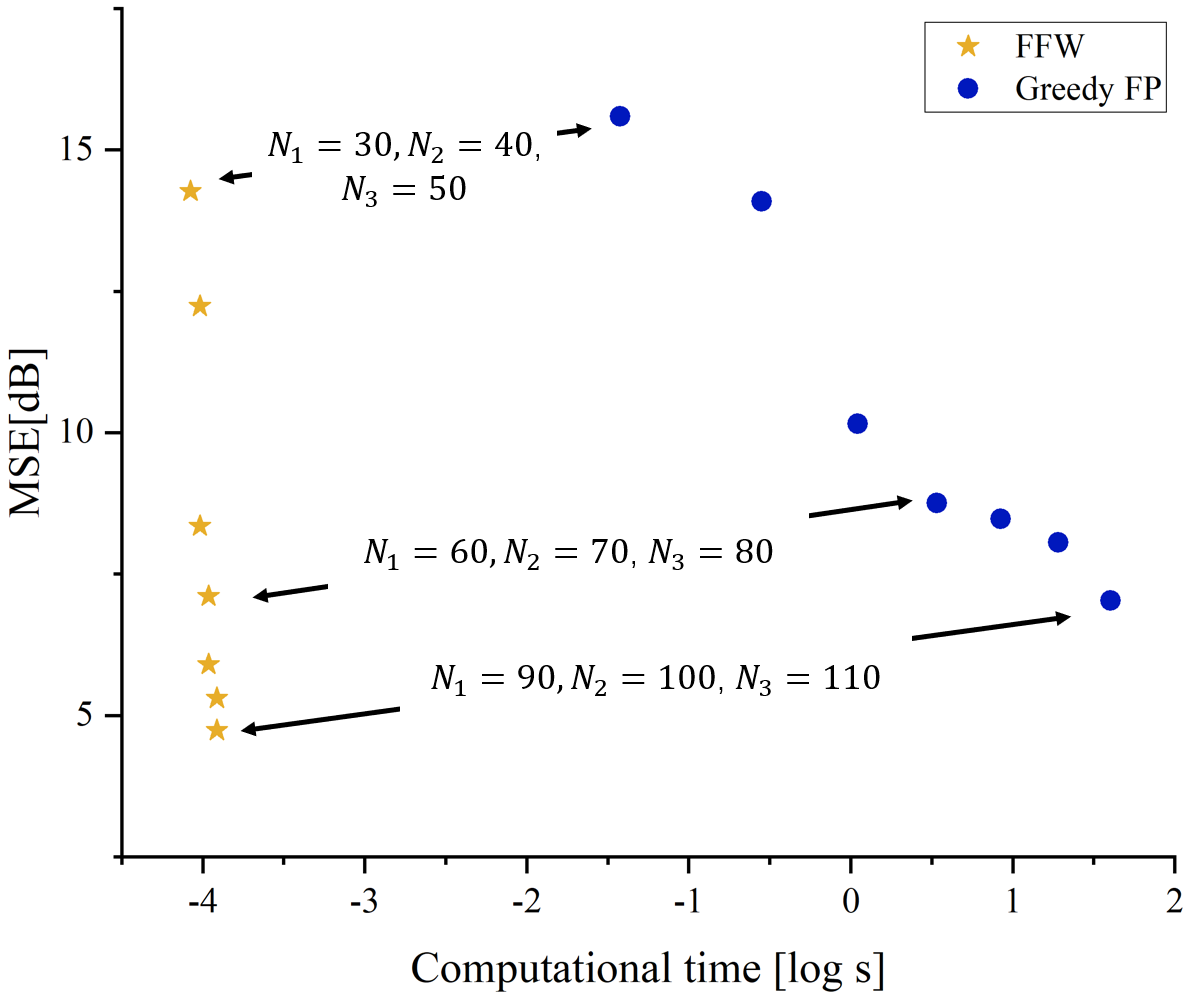}%
    \label{2_b}}
  \caption{Analysis of the tradeoff between computational time and MSE for FFW with FMBS, FrameSense and  Greedy FP}
  \label{fig2}
\end{figure*}

In the fourth experiment, we randomly draw six matrices $\mathbf{\Psi }\in \mathbb{R} ^{N\times K}$ with $K=40$ and $N \in \{100,150,200,250,300,350,400\}$, while we place $L = 0.5N$ sensors for vector signals. For tensor signals, we draw  matrices $\left\{ \mathbf{\Psi }_i\in \mathbb{R} ^{N_i\times K_i} \right\} _{i=1}^{R=3}$ with dimensions $K_1=10$, $K_2=20$ and $K_3=15$, while $
  \left( N_1,N_2,N_3 \right) \in \left\{ \left( 30+\omega ,40+\omega ,50+\omega \right) |\omega =0,10,20,30,40,50,60 \right\}
$.
We measure the  computational time together with the  MSE, showing that for vector signals, while FFW is significantly faster than FMBS and  FrameSense, the difference in MSE is minimal (Figure \ref{2_a}). Moreover, when sampling tensor signals, the gap of  computational time between FFW and Greedy FP is greater, and the MSE performance of FFW is also better than Greedy FP (Figure \ref{2_b}).

\begin{figure*}[!t]
  \centering
  \subfloat[Vector signals.]{\includegraphics[width=2.5in]{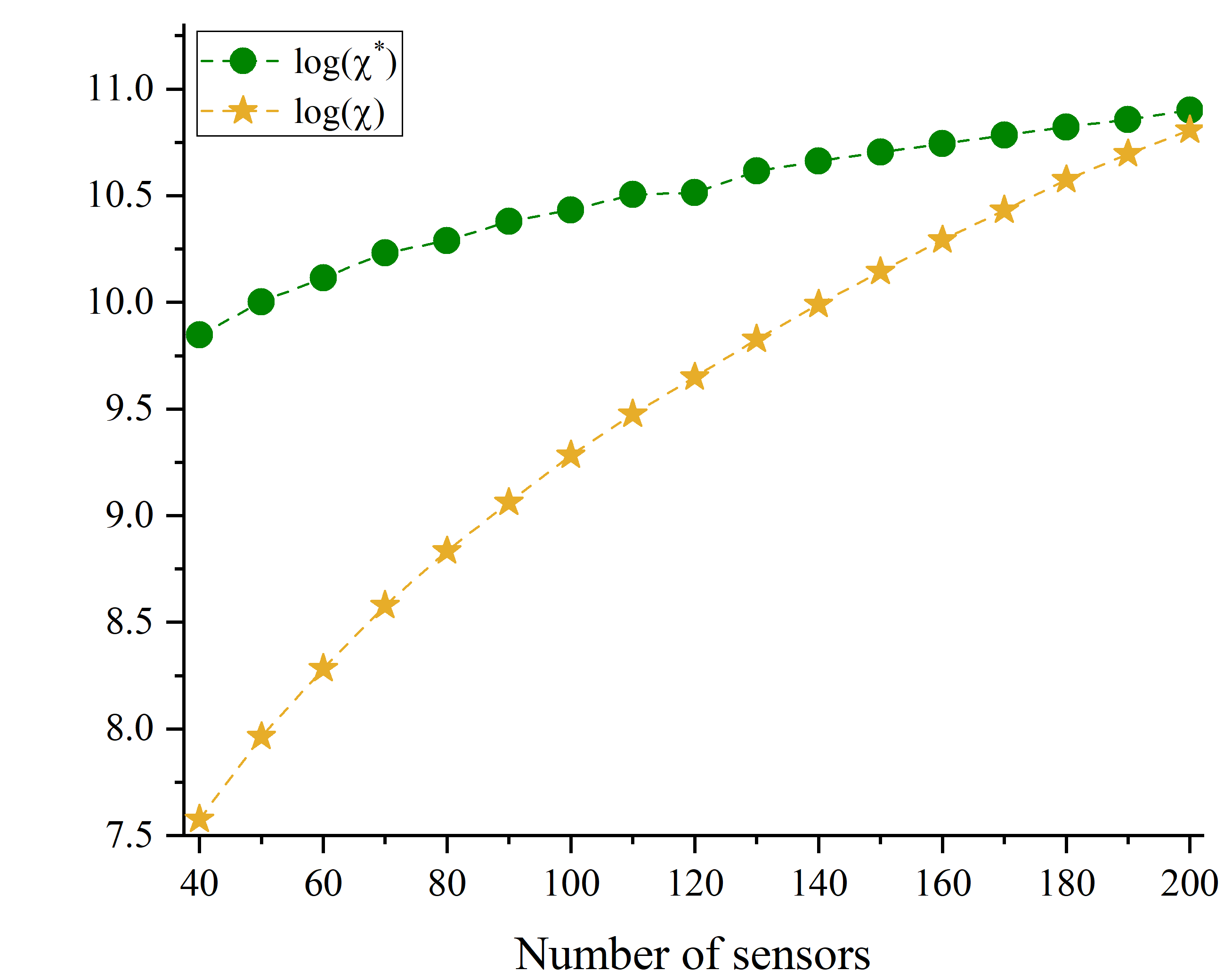}%
    \label{6_a}}
  \hfil
  \subfloat[Tensor signals.]{\includegraphics[width=2.5in]{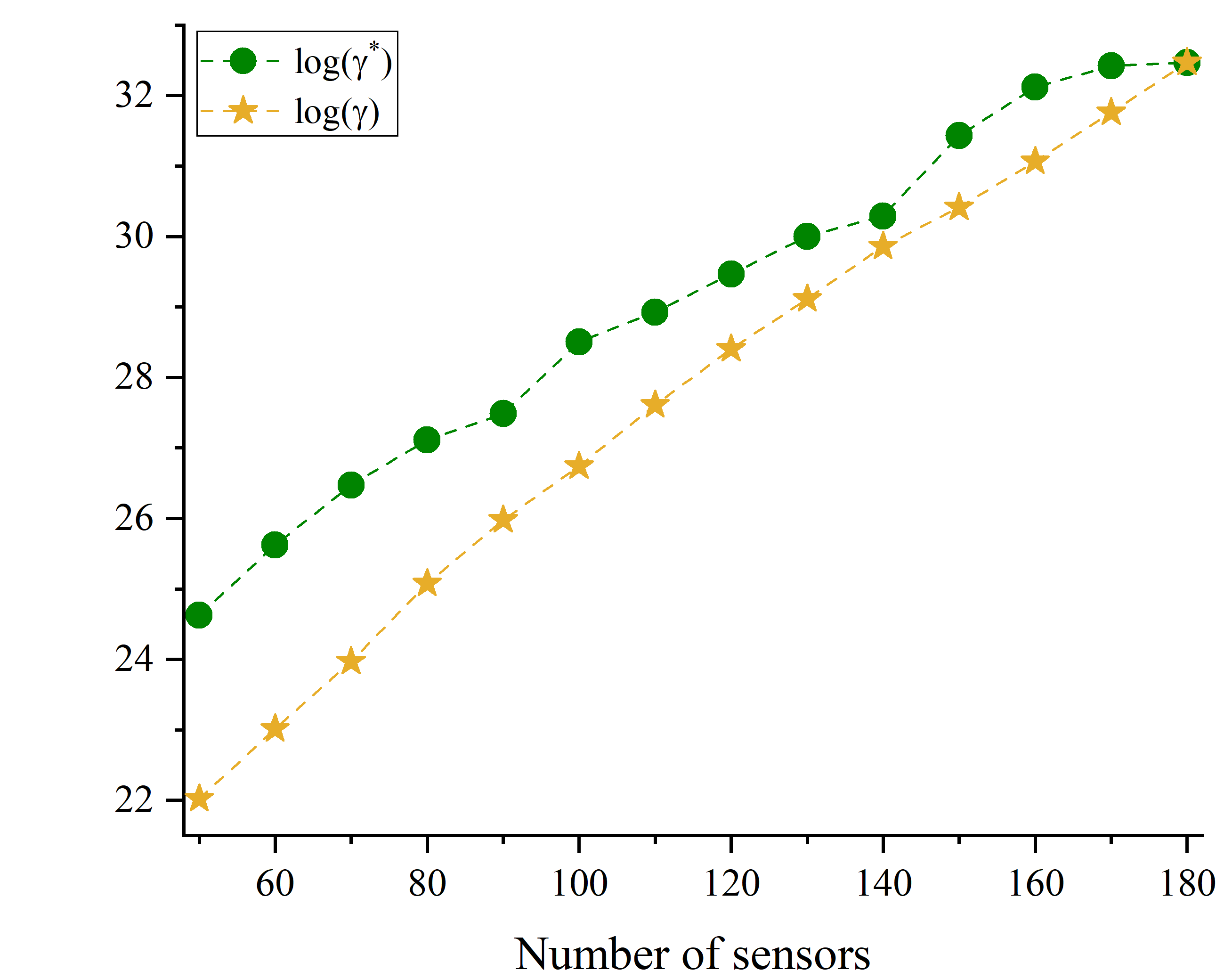}%
    \label{6_b}}
  \caption{The comparison between the bound  given by Theorem \ref{theo2} (Theorem \ref{theo3}) and $F\left( \mathcal{L} ^{\prime} \right)$, where $\mathcal{L} ^{\prime}$ is the solution obtained by Algorithm \ref{alg1} (Algorithm \ref{alg2}).  }
  \label{fig6}
\end{figure*}

In the final experiment of this subsection, we verified the approximation factors provided by Theorem \ref{theo2} and Theorem \ref{theo3}. Firstly, we used the data generated in the first experiment to validate Theorem \ref{theo2}. Specifically, let $\chi$ and $\chi^*$ be $F\left( \mathcal{L} ^{\prime} \right)$ and $\gamma F\left( \mathcal{L} ^* \right)
$ in Theorem \ref{theo2}.   Since $\mathcal{L} ^*$ can not be obtained  in practice, we use the optimal result of random sampling 10000 times to instead of $\mathcal{L} ^*$. The experimental results are shown in Figure \ref{6_a}
Then, we employed the data generated in  the second experiment  to verify Theorem \ref{theo3}.  Let $\gamma $ and $\gamma ^*$ be $F\left( \mathcal{L} ^{\prime} \right)$ and  $e^{M-\frac{M^2}{2R}}F\left( \mathcal{L} ^* \right)$ in Theorem \ref{theo3}, and the experimental results are also depicted in Figure \ref{6_b}. The yellow points in these figures represent $\log (\chi)$($\log (\gamma)$), and the green represent $\log (\chi^*)$($\log (\gamma^*)$).
As the value of $L$ grows, we observe a trend that $\chi$($\gamma $) becomes closer to $\chi^*$($\gamma^* $). Additionally, when $L$ equals $N$, it is observed that $\gamma$ becomes equivalent to $\gamma^* $.

\subsection{Sampling and reconstruction of image data}

In this subsection, we used a classical demo image in MATLAB  called “Peppers” for performance comparison. We first used  FFW  and Greedy FP for sampling real-world image data and use the least squares method to reconstruction. Then, we randomly select $L$ rows and columns of image data,  and apply Convolutional Conditional Neural Process (ConvCNP) \cite{CNP} as the model to completion.

\begin{figure*}[!t]
  \centering
  \subfloat[Original image]{\includegraphics[width=2in]{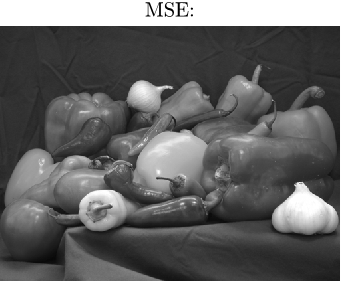}%
    \label{3_a}}
  \hfil
  \subfloat[FFW]{\includegraphics[width=2in]{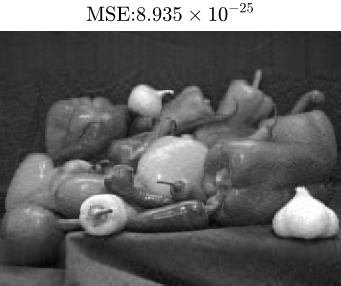}%
    \label{3_b}}
  \hfil
  \subfloat[Greedy FP]{\includegraphics[width=2in]{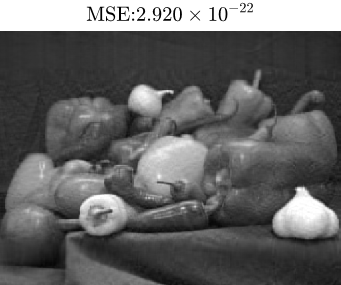}%
    \label{3_c}}
  \caption{Visualization of different sampling and reconstruction results of real-world image (“Peppers”).}
  \label{fig3}
\end{figure*}

Firstly, we decompose the image $\mathbf{X}$   into $\mathbf{U}_1 \mathbf{G} \mathbf{U}_2^{\top}$, where $\mathbf{U}_1$  and $\mathbf{U}_2^{\top}$   are non-orthogonal matrices. Then we use the first $K_1$ columns of $\mathbf{U}_1$ and $K_2$ column of $\mathbf{U}_2$ as two factor matrices for sampling and reconstruction. We set the dimension of low-dimensional parameter matrix to be $40 \times 40$, i.e., $K_1=40$ and $K_2=40$. Let the number of sensors $L=400$. The results are  showed in Figure \ref{fig3}. One can see that MSE of sampling by Greedy FP is much bigger than MSE of sampling by FFW.

Then, we randomly sample the image data, and      complete it by ConvCNP. The number of sensors $L$ is $400, 600, 700 $ and $800$ respectively.

ConvCNP combines a convolutional neural network with features of a Gaussian process to model the conditional distribution of inputs as a way of dealing with uncertainty modeling between input-output pairs. Compared to traditional point-based prediction models, ConvCNP is able to make high-quality predictions on unobserved data by modeling conditional distributions with a small amount of observed data. ConvCNP is suitable for a range of canonical machine learning tasks, including regression, classification and image completion \cite{CNP}.

For image completion task, we choose a commonly used dataset -— the CIFAR10 dataset \cite{zhou22}, which is  a commonly used dataset for image completion task. The CIFAR10 dataset is of moderate size, containing 10 categories of color images with 6000 $32 \times 32$ RGB images in each category. These categories cover different items and scenes in daily life, which can enable the generative model to learn the features among different categories and generate images with generalization. Using the public data set CIFAR10, we train the ConvCNP model for 200 epochs using Adam \cite{zhou23} with learning rate of $5 \times 10^{-4}$, batch size of 256. We input the original image “Peppers” and control the proportion of the original image to be masked by setting the masking factor. Since the model needs to complement the images masked by random whole rows and columns in the testing phase, we train the model by randomly selecting whole rows and columns to mask images. After the training is completed, we save the optimal model parameters. Then we convert the masked matrix into a grayscale image to match the number of channels (the number of channels in CIFAR10 dataset is 3), and then input the image into the model as an observation for image completion. The results can be seen in Figure \ref{fig4}.


Obviously, the results completed by convCNP is not as good as those reconstructed by least squares method (Figure \ref{fig3}). Usually, the predictive performance of an intelligent model is highly dependent on the number of context points \cite{zhou24,zhou25,zhou26}.
While the ConvCNP model effectively predicts masked regions with a few context points, the uncertainty in its predictions increases as the number of context points decreases.
In the experiments, the model omits a lot of contextual information by masking images with whole rows and columns, resulting in reducing the prediction accuracy of the model. In addition, intelligent models require high computational resources and high time costs. However, under the model of this paper, using the singular value decomposition of the image data, the core is reconstructed by the least square method after fast sampling, so that the original image can be reconstructed more quickly and accurately.

\begin{figure*}[!t]
  \centering
  \subfloat{\includegraphics[width=1.5in]{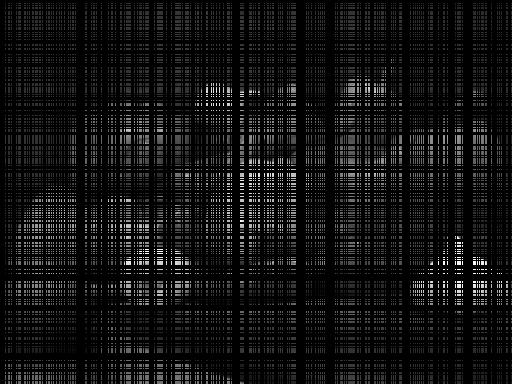}}
  \hfil
  \subfloat{\includegraphics[width=1.5in]{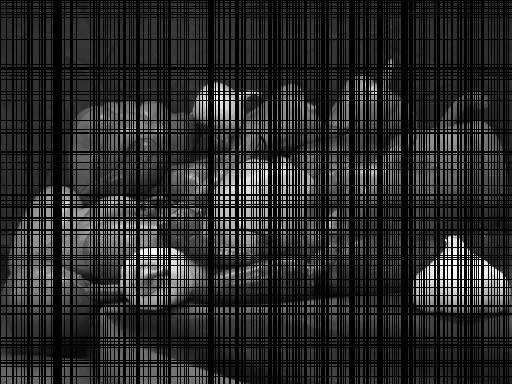}}
  \hfil
  \subfloat{\includegraphics[width=1.5in]{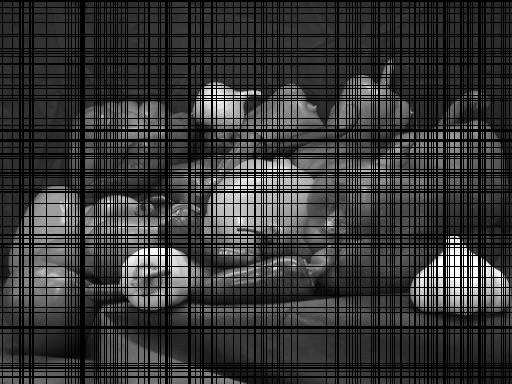}}
  \hfil
  \subfloat{\includegraphics[width=1.5in]{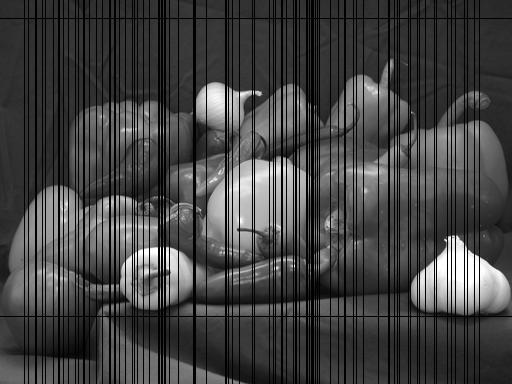}}
\end{figure*}
\begin{figure*}[!t]
  \centering
  \subfloat[$L=400$]{\includegraphics[width=1.5in]{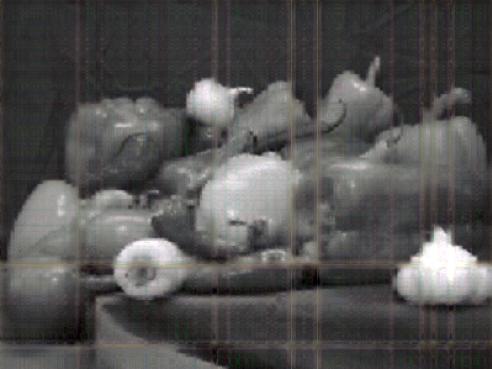}%
    \label{4_a}}
  \hfil
  \subfloat[$L=600$]{\includegraphics[width=1.5in]{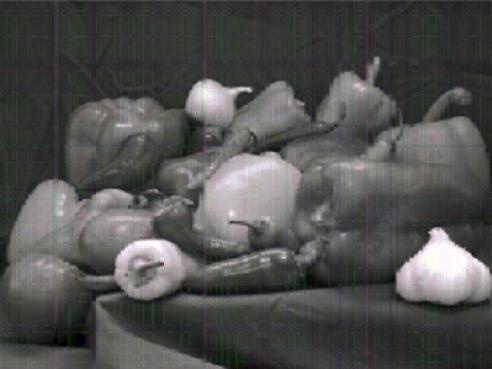}%
    \label{4_b}}
  \hfil
  \subfloat[$L=700$]{\includegraphics[width=1.5in]{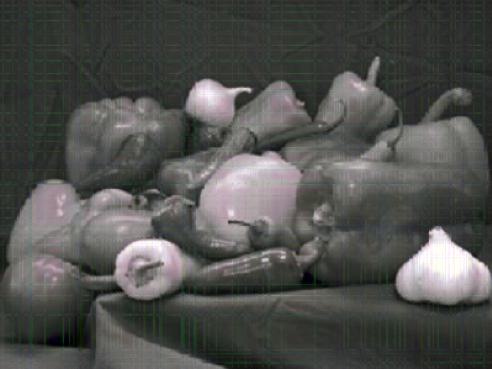}%
    \label{4_c}}
  \hfil
  \subfloat[$L=800$]{\includegraphics[width=1.5in]{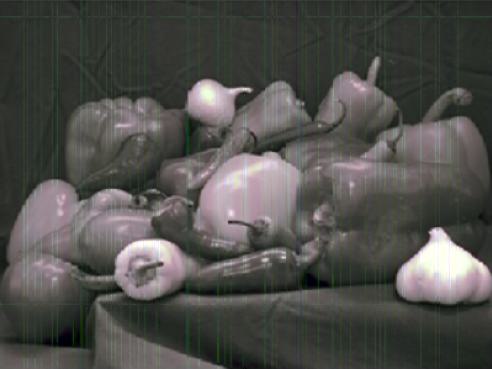}%
    \label{4_d}}
  \caption{Visualization   of using convCNP completion after random sampling of real-world image (“Peppers”).}
  \label{fig4}
\end{figure*}
\section{Conclusion}

In this paper, we proposed a fast algorithm (FFW) to sampling for inverse problems with vectors and tensors.  We first provide the closed-form expressions of multilinear extensions for the frame potential of  pruned  matrices (Theorem \ref{theo0}). For faster sampling, we  design  FFW to sampling  vectors  with  complexity $\mathcal{O}(NK^2)$, and extend FFW to sampling high-order tensors  with complexity $\mathcal{O} \left( N_{\max}K_{\max}^2 \right)$. We also give the approximation factor of FFW for a special class of factor matrices (Theorem  \ref{theo2} and Theorem  \ref{theo3}). Then, we experimentally demonstrate the strength of FFW in performance and running time, compared with FrameSense, FMBS and Greedy FP, and verify
that FFW also has higher performance under the cases that factor matrices do not satisfy the conditions of Theorem \ref{theo1} and Theorem \ref{theo3}.
Finally, for image data, the results of using FFW sampling  and least squares reconstruction are better than the results of using convCNP completion after random sampling.

\bibliography{bare_jrnl_new_sample4}
\bibliographystyle{IEEEtran}

\end{document}